\documentclass[a4paper]{amsart}
\usepackage{graphicx}
\usepackage{amssymb}
\usepackage{epstopdf}
\input xy
\xyoption{all}

\newtheorem{lemma}{Lemma}[section] 
\newtheorem{propos}[lemma]{Proposition}

\newtheorem{theorem}[lemma]{Theorem}
\newtheorem{cor}[lemma]{Corollary}

\newcommand{\CM}{\hbox{{$\mathcal M$}}}

\newcommand{\CO}{\hbox{{$\mathcal O$}}}

\newcommand{\CR}{\hbox{{$\mathcal R$}}}
\newcommand{\cg}{\hbox{{$\mathfrak g$}}}

\newcommand{\C}{\mathbb{C}}
\newcommand{\R}{\mathbb{R}}

\newcommand{\Z}{\mathbb{Z}}

\newcommand{\del}{\partial}

\newcommand{\extd}{\mathrm{d}}

\newcommand{\isom}{{\cong}}
\newcommand{\eps}{{\epsilon}}
\newcommand{\tens}{\mathop{\otimes}}
\newcommand{\la}{{\triangleright}}
\newcommand{\ra}{{\triangleleft}}

\newcommand{\id}{\mathrm{id}}

\newcommand{\<}{\langle}
\renewcommand{\>}{\rangle}

\renewcommand{\o}{{}_{\scriptscriptstyle(1)}}
\renewcommand{\t}{{}_{\scriptscriptstyle(2)}}

\begin{document}
\title{Nonassociative Riemannian geometry by twisting}
\author{Edwin Beggs \& Shahn Majid}
\date{8/12/09}                                           
\address{EJB: Department of Mathematics, Swansea University\\
Singleton parc, Swansea SA2 8PP, UK}
\address{SM: School of Mathematical Sciences, Queen Mary University of London\\
327 Mile End Rd, London E1 4NS\\ }
\email{e.j.beggs@swansea.ac.uk \& s.majid@qmul.ac.uk}
\begin{abstract}Many quantum groups  and quantum spaces of interest can be obtained by cochain (but not cocycle) twist from their corresponding classical object. This failure of the cocycle condition implies a hidden nonassociativity in the noncommutative geometry already known to be visible at the level of  differential forms. We extend the cochain twist framework to connections and Riemannian structures and provide examples including twist of the $S^7$ coordinate algebra to a nonassociative hyperbolic geometry in the same category as that of the octonions. \end{abstract}
\maketitle
\section{Introduction}

The idea of noncommutative geometry, of course, is to replace points in an actual space by a coordinate algebra. These are generally generated by variables $x,y,z$ ... say  enjoying algebraic properties for  addition, multiplication and other relations paralleling the numbers for which they could be viewed as placeholders. Thus, solving $x^2+y^2+y^3=1$ over $\R$ would bring one to the actual points of a sphere, but one could also consider this equation more abstractly with $x,y,z$ as generators of an algebra.  And if in place of the usual $xy=yx$ etc we have some other non-commutation relations then, clearly, these generators could never be realised as numbers. Over the last two decades it has become obvious that the assumption of commutativity is a historical accident to do with the fact that classical mechanics was discovered before quantum mechanics and there is no particular reason, certainly from a mathematical perspective, not to allow noncommutavity as a generalisation of usual geometry more applicable to the quantum world. Physicists are also somewhat familiar with noncommuting variables from quantum mechanics and there are  physical reasons to think that  plausibly the first quantum gravity corrections to geometry should be expressed as such noncommutativity\cite{MaSch}. Accordingly, a great deal of effort was been put into developing such a geometry both in Connes `spectral triple' approach and in a `quantum groups' approach modelled around quantum groups such as the deformations $C_q(G)$ and their homogeneous spaces in the first instance. Both approaches are by now somewhat mature and the quantum groups approach includes specific testable predictions for Planck scale physics\cite{AmeMa}. Meanwhile, the spectral triples approach includes specific  predictions and mass relations for the standard model\cite{ConCham}. 

In this article we go one step further and consider geometry that is nonassociative in the same algebraic sense. While noncommutative geometry is motivated from quantum theory and hence somewhat familiar, this is not so for nonassociativity. If we simply drop associativity then most of what we do in algebra and geometry becomes vastly more complicated if not impossible. Therefore the first point is that we should only approach nonassociative geometry under duress. Indeed, common to many approaches to noncommutative geometry is the idea of replacing differential structures by an algebraic formulation of `exterior algebra of differential forms'  with generators $\extd x,\extd y,\extd z$... say,  formalized as an associative  differential graded algebra $(\Omega,\extd)$ over our possibly noncommutative  `coordinate algebra'. This is pretty much a prerequisite to do any kind of physical model building on noncommutative geometry, wave equations, maxwells equations etc. Unfortunately, after about 20 years experience with model building it turned out to be {\em not quite} possible to do this while maintaining a strict correspondence with classical differentials and classical symmetries. When quantization fails in this respect one speaks of a quanutm anomaly, and  in previous work we have similarly noted a fundamental anomaly or obstruction in the entire programme of noncommutative geometry.

\begin{theorem} \cite{BM1} The standard quantum groups $C_q(G)$ of simple Lie algebras $\cg$ do {\em not} admit an associative differential calculus $\Omega(C_q(G))$ in deformation theory with left and right translation-covariantce and classical dimensions.
\end{theorem}

We showed that this obstruction, like other anomalies in physics, can be expressed at the semiclassical level as curvature, in our case of a certain Poisson preconnection, and that for $\cg$ semisimple there is no flat connection of the required type. The same applies to enveloping algebras $U(\cg)$ when $\cg$ is semisimple, viewed as quantisation of $\cg^*$ with its Kirillov-Kostant Poisson bracket:

\begin{theorem} \cite{BM2} The classical enveloping algebras $U(\cg)$ of  simple Lie algebras $\cg$ do {\em not} admit a differential calculus $\Omega(U(\cg))$ in deformation theory with ${\rm ad}$-invariance under $\cg$ and  classical dimensions.
\end{theorem}

This was again proven using similar methods, notably Kostant's invariant theory. We note that an important first attempt at such a rigidity result was in \cite{Haw} which, in particular, pointed to the role of preconnections. Moreover, attempts at such calculi at the quantum group level had been found to require extra dimensions in the cotangent bundle and this was understood now as a way to `neutralise' the anomaly, again much as for other anomalies in physics. Here \cite{BM1,BM2} also provided an alternative: using cochain twisting methods  it was showed that we can always keep classical dimensions and deform the classical picture, provided we allow {\em nonassociative geometry}. In short, if one wants a strict deformation-theoretic correspondence on quantization, we have
\[ NCG \Rightarrow NAG\]
(nonassociative geometry implies nonassociative geometry). In the present paper we explore this radical alternative further. And if we want to be speculative then, once we allow nonassociative exterior algebras, we should also alllow nonassociative coordinate algebras, i.e. spaces themselves not only their differential geometry could be nonassociative. 

In entering this nonassociative world we also need tight controls. The key idea is to work in a `quasiassociative' setting in which coordinate algebras are nonassociative but in a tightly controlled way by a {\em multiplicative associator}. In mathematical terms this means that the algebra {\em is} associative but in some monoidal category different from that of vector spaces. A theorem of Maclane says that all constructions in such a category can be done {\em as if} associative i.e. without worrying about brackets. One must then insert the associator $\Phi_{X,Y,Z}:X\tens (Y\tens Z) \to (X\tens Y)\tens Z$ between any three objects as needed in order to make sense of expressions and Maclane's theorem says that in a monoidal category any different ways to insert $\Phi$ as needed will give the same result.   We focus in this article on a specific `twist quantisation functor' that results in noncommutative algebras in such a monoidal category and also quantises all other covariant structures with reference to a classical symmetry.   

To first explain the background,   V.G. Drinfeld \cite{Dri:qua} showed that all quantum group enveloping algebras $U_q(\cg)$ could be obtained as follows: start with the classical symmetry algebra $U(\cg)$ {\em but} viewed as a quasi-Hopf algebra with respect to a certain element $\phi\in U(\cg)^{\bar{\tens} 3}$ (a topological tensor product) obtained by solving the KZ equations. Its category of representations is a braided monoidal category with $\Phi$ giveny the action of $\phi$ (and braiding given by the action of $\CR=e^{h t}$ where $q=e^{h\over 2}$ and $t$ is the split casimir of $\cg$). Drinfeld showed that there exists a certain cochain $F\in U(\cg)^{\bar{\tens} 2}$  such that conjugating the classical coproduct of $U(\cg)$ by $F$ gives $U_q(\cg)$. Here $\phi$ also twists to some $\phi_F$ but $F$ is chosen so that $\phi_F=1$ -- an ordinary Hopf algebra. In \cite{AkrMa,BM1} we developed a rather different application of the same data, namely start with a completely classical data $(U(\cg),\phi=1,\CR=1)$ and twist this by Drinfelds $F$ (conjugate the classical coproduct). The result $U(\cg)^F$ looks like $U_q(\cg)$ {\em but} regarded as a quasi-Hopf algebra with $\phi_F$ closely related to Drinfelds and $\CR_F=F_{21}F^{-1}$ making its category of representations a (symmetric) monoidal category.     There is also an important spin-off\cite{Ma:book}: if the classical object $U(\cg)$ acts on a classical coordinate $C(M)$ in algebraic terms then because of the functorial nature of the twist construction $C(M)$ also gets `quantized' to $C(M)_F$, typically noncommutative and nonassociative because it lives in the category of representation of $U(\cg)^F$. Viewed in that category it is in fact commutative and associative, but the category is no longer the usual one due to nontrivial $\Phi$. Similarly the classical $\Omega(M)$ gets quantised to $\Omega(M)_F$ as a calculus on $C(M)_F$. In this way any classical geometry can be systematically quantised with respect to a classical symmetry and a choice of cochain $F$ -- provided we can live with potential nonassociativity.  

Note that in geometry we are interested not in enveloping algebras but in coordinate algebras. Therefore, while not essential, we  will convert the above to  work with classical $C(G)$ and with $C_q(G)$ viewed as a coquasiHopf algebra $C(G)^F$ and  replace the action of $U(\cg)$ on $C(M)$ by a {\em coaction} of $C(G)$. In particular $G$ acts ($C(G)$ coacts) on $C(G)$ by left translation and induces a cochain quantisation $C(G)_F$ covariant under $C(G)^F$\cite{AkrMa,BM2}. We recall this less familiar setting in the Preliminaries below. Depending on the choice of cochain some of the twisted algebras may remian accidentally associative -- this would be the case for Drinfelds cochain and $C(G)^F$ -- and in that case the nonassociativity is hidden. But it is still present and appears typically in  the differential calculus even in this case. This programme has already been carried quite far and covers the general principles and the calculus, while in the present article we now study how the next main layers of geometry, connections and curvature etc, behave under such cochain twists. The short answer is that everything works as expected if everything is covarant under the twisting classical symmetry. This also points to the limitations of the twist approach: it works {\em too well} -- different but covariant ways to express classical constructions will all twist and only by having a general theory of noncommutative (and now nonassociative) geometry that makes sense {\em beyond} examples given by twisting will we be able to know which of these is most natural. We will illustrate this with the Ricci tensor in Section~7. One could also consider metrics, for example, that are  not invariant under the classical symmetry, although we do not do so here.

Finally, for nonassociative geometry we need concrete examples to build up our experience. About the only example that most readers may have some experience with is the octonions  and these provided the first concrete example of the above `cochain quantisation' in action \cite{AlbMa}. The octonions are a dimension 8 `algebra' but with a nonassociative product, which is, however, alternative, in the sense that if two elements are repeated then they associate, so $x(yx)=(xy)x$ etc. The full picture coming out of the quantum group theory amounts to the following. One can choose a  basis $\{e_a\}$ labelled by $a\in \Z_2^3$ (a 3-vector with values in $\Z_2$) such that\cite{AlbMa}
\[ (e_a e_b)e_c=\phi(a,b,c) e_a(e_b e_c),\quad \phi(a,b,c)=(-1)^{a\cdot (b\times c)}=(-1)^{|a b c|}.\]
This means that three elements associate {if and only if} their labels are linearly dependent as vectors over $\Z_2$, otherwise there is a -1. This was obtained in \cite{AlbMa} precisely as a cochain quantization which we can view as $C(\Z_2^3)_F$ of $C(\Z_2^3)$ under the action of $\Z_2^3$ by left translation and a certain choice of $F$. One can then proceed to differential geometry on the octonions as coordinate algebra\cite{Ma:gau}.   Our reinterpretation of Drinfeld's ideas to produce $C(G)_F$ is inspired by this and exactly parallels the construction of the octonions but now for a Lie group. The octonion example also illustrates that even though twisting is `easy' in the sense that everything is somewhat automatic by an equivalence of categories reflected in the fact that $\phi$ is cohomologically trivial, this does not mean that the resulting objects are uninteresting. The octonions, for example, are a division algebra while the group algebra of $\Z_2^3$ is certainly not. At the same time, once we adopt our framework of geometry in monoidal categories, we are not limited to such twist or coboundary examples but can consider more general ones. We are also not  limited to differential geometry and can consider nonassociative quantum mechanics. Some recent suggestions about the latter are in \cite{Jens}. It would also be interesting to apply non-associative spectral triples to the standard model, for example based now on octonions rather than two copies of the quaternions for the `internal geometry'.

\subsubsection*{Note added.} Some of Tony Sudbery's work was towards the noncommuative geometry coming out of quantum groups and it is therefore a pleasure to contribute the present article to his Festschrift volume. There is also a degree of irony involved. In the mid 1990s Tony  left quantum groups to work on the seemingly unrelated area of octonions and other topics, but quantum groups somehow followed him with results on octonions as above. Meanwhile, those who stayed behind to work on purely on quantum groups and their noncommutative geometry are now being forced  to learn lessons on how to do it from  the octonions.

\section{Preliminaries: cochain twists}

We work over a field $k$ but for physics we will often specialise this to $\C$ or $\R$. We will need to use quantum group or Hopf algebra methods and we refer to \cite{Ma:book} for an introduction. We recall that a quantum group $H$ is a bialgebra equipped with a notion of `linearised inverse' or `antipode' $S:H\to H$. Here a bialgebra means that the algebra $H$ is also a coalgebra in a compatible way. A concept of a coalgebra is just the same as the concept of an algebra but with arrows reversed, so there is a `coproduct' $\Delta:H\to H\tens H$ and `counit' $H\to k$, and we require $\Delta$ (and hence $\eps$) to be algebra homomorphisms. There is a standard (Sweedler) notation $\Delta h=h\o\tens h\t$ which we shall use  for the result of $\Delta$ as an element of $H\tens H$ (summation understood).

Suppose that $H$ is a coalgebra. We consider its category of left comodules
${}^H\mathcal{M}$ and we write the coaction on an object $V$
as $\lambda:V\to H\tens V$, or as $v\mapsto v_{[-1]}\tens v_{[0]}$ with summation understood. The concept of a comodule here is just the same as that of a module but with arrows reversed. 
If  $H$ is a Hopf algebra or we can use its algebra 
product to define a tensor product on ${}^H\mathcal{M}$ by
\begin{eqnarray*}
\lambda(v\tens w)\,=\, v_{[-1]}\,w_{[-1]}\tens (v_{[0]}\tens w_{[0]})\ .
\end{eqnarray*}
The identity in the category is just the underlying field $k$, with $\lambda:k\to H\tens k$
being $1_H\tens$. 
Note that at this stage, a bialgebra structure on $H$ would suffice, but it shall be convenient to have a Hopf algebra later. An algebra in the category ${}^H\mathcal{M}$ is an object $A$ with
an associative product
$\mu:A\tens A\to A$ and a unit $1_A:k\to A$ which are $H$-comodule maps. 

If $G$ is a finite or algebraic group we are interested in the point of view where $H$ could be its function or coordinate algebra $k(G)$ with its usual pointwise product. If similarly $A=k(X)$ is the algebra of functions one some set or algebraic variety $X$ and $G$ acts on $X$, this would be expressed in our algebraic terms as a coaction $\lambda:A\to H\tens A$. For example, in the finite group case
 \begin{eqnarray*}
\lambda(f)(x)\,=\,\sum_{g\in G} \delta_g.f(gx)\ 
\end{eqnarray*}
where $\delta_g$ is the Kronecker $\delta$-function. It may seem strange to do everything `backwards' in terms of coactions but this will allow us to work with algebraic tensor products and algebraic groups. Our framework also allows  $A$ and $H$ to be noncommutative in which case there would be no actual spaces $X,G$ respectively.

Now we shall consider a different tensor product $\tens^F$ for ${}^H\mathcal{M}$. 
We shall do this by changing the product on $H$ while leaving the coproduct unaltered. 
As the coproduct is the same, we can consider the category ${}^H\mathcal{M}$
to be the same, we have only changed the tensor product. We begin with a 
linear map $F:H\tens H\to k$ which is invertible, in the sense that there
is another map $F^{-}:H\tens H\to k$ so that
\begin{eqnarray*}
F(a_{(1)}\tens b_{(1)})\, F^{-}(a_{(2)}\tens b_{(2)})\, =\,
F^{-}(a_{(1)}\tens b_{(1)})\, F(a_{(2)}\tens b_{(2)})\, =\, \epsilon(a)\,\epsilon(b)\  
\end{eqnarray*}
and which obeys $F(a\tens 1)=F(1\tens a)=\eps(a)$. Such a map is called a 2-{\em cochain} on $H$\cite{Ma:book}. Now we specify the tensor product $\tens^F$ on ${}^H\mathcal{M}$ by  the same vector space $V\tens W$ as usual but with coaction
\begin{eqnarray*}
\lambda(v\tens^F w)\,=\, v_{[-1]}\circ w_{[-1]}\tens (v_{[0]}\tens^F w_{[0]})\ ,
\end{eqnarray*}
where the modified product $\circ:H\tens H\to H$ is given in terms of the original product by
\begin{eqnarray*}
a\circ b \,=\, F(a_{(1)}\tens b_{(1)})\, a_{(2)} b_{(2)}\, F^{-}(a_{(3)}\tens b_{(3)})\, .
\end{eqnarray*}
This product may not be associative -- in general $H^F=(H,\circ)$ is  coquasi-Hopf algebra\cite{Ma:book}, but this is not a problem in defining $\tens^F$. 
Where it does appear as a complication is that the trivial `identity' map from
$(U\tens^F V)\tens^F W\to U\tens^F (V\tens^F W)$ is not a morphism.
Instead we can use the non-trivial associator
$\Phi_{U,V,W}:(U\tens^F V)\tens^F W\to U\tens^F (V\tens^F W)$ defined by
\begin{eqnarray*}
\Phi((u\tens^F v)\tens^F w) &=& (\partial F)^-(u_{[-1]}\tens v_{[-1]}\tens w_{[-1]}) \ 
u_{[0]}\tens^F (v_{[0]}\tens^F w_{[0]})\ ,\cr
\Phi^{-1}(u\tens^F (v\tens^F w)) &=& (\partial F)(u_{[-1]}\tens v_{[-1]}\tens w_{[-1]}) \ 
(u_{[0]}\tens^F v_{[0]})\tens^F w_{[0]}\ ,
\end{eqnarray*}
where $\partial F,(\partial F)^-:H\tens H\tens H\to k$ is defined by
\begin{eqnarray*}
(\partial F)^-(a\tens b\tens c) &=& F(a_{(1)}\tens b_{(1)})\, F(a_{(2)}b_{(2)}\tens c_{(1)})\, 
F^{-}(a_{(3)}\tens b_{(3)} c_{(2)})\,
F^{-}(b_{(4)}\tens c_{(3)})\ ,\cr
(\partial F)(a\tens b\tens c) &=& 
F(b_{(1)}\tens c_{(1)})\,
F(a_{(1)}\tens b_{(2)} c_{(2)})\, 
F^-(a_{(2)}b_{(3)}\tens c_{(3)})\, 
F^-(a_{(3)}\tens b_{(4)})\ .
\end{eqnarray*}
It would be customary here to restrict to the case where $F$ is a cocycle, i.e.\
$(\partial F)(a\tens b\tens c)=\epsilon(a)\epsilon(b)\epsilon(c)$. This
means that the new tensor product would also be trivially associated. 
However it is not really necessary to suppose this, if the reader does not mind dealing with nontrivially associated tensor categories, and we shall not do so. 

It will be useful to note that there is a natural transformation $c$ between
the two tensor product functors, $\tens$ and $\tens^F$, from ${}^H\mathcal{M}
\times {}^H\mathcal{M}\to {}^H\mathcal{M}$, given by
\begin{eqnarray} \label{kugdscvahksuv}
c(v\tens w) &=& F^-(v_{[-1]}\tens w_{[-1]})\, v_{[0]}\tens^F w_{[0]}\ 
\end{eqnarray}
which, together with the identity map on objects and morphisms provids an equivalence of monoidal categories $({}^H\CM,\tens)$ and $({}^H\CM,\tens^F)$. The latter is the category of comodules of $H^F=(H,\circ)$.

\section{Twisting algebras and modules}

An algebra $(A,\mu)$ in the tensor category $({}^H\mathcal{M},\tens)$ is twisted to
an algebra $A_F=(A,\mu^F)$ in $({}^H\mathcal{M},\tens^F)$ by
\begin{eqnarray}\label{bullet}
\mu^F(a\tens b)\,=\, a\bullet b\,=\, F(a_{[-1]}\tens b_{[-1]})\ a_{[0]}b_{[0]}\ .
\end{eqnarray}
A quick check will show that this is associative in $({}^H\mathcal{M},\tens^F)$,
remembering to change the order of the bracketing using $\Phi$. Thus,
\begin{equation}\label{bulletassoc} (a\bullet b)\bullet c= (\partial F)^-(a_{[-1]}\tens b_{[-1]}\tens c_{[-1]}) \ 
a_{[0]}\bullet(b_{[0]}\bullet c_{[0]}).\end{equation}

\medskip
 We can also
twist left or right $A$-modules $V$ in $({}^H\mathcal{M},\tens)$
(meaning that the left or right actions have to preserve the $H$-coaction) to get left or right modules of $(A,\mu^F)$ 
\begin{eqnarray*}
a\la^F v\,=\, F(a_{[-1]}\tens v_{[-1]})\ a_{[0]}\la v_{[0]}\ ,\quad
v\ra^F a\,=\, F(v_{[-1]}\tens a_{[-1]})\ v_{[0]}\ra a_{[0]}\ .
\end{eqnarray*}
for $v\in V$ and $a\in A$, and where we use $\la$ and $\ra$ for the left and right actions respectively. These obey the requirements of an action of $(A,\mu^F)$ in the category with $\tens^F$, i.e. inserting $\Phi$ in the form $\del F^-$ much as in (\ref{bulletassoc}). When we have both left and right actions on $V$ we can require them to commute so $(a\la v)\ra b=a\la (v\ra b)$. We say that $V\in {}_A\CM_A$, the category of bimodules. If in addition our bimodule is left covariant under $H$ then twisting both actions gives us a bimodule of $(A,\mu^F)$ with respect to $\tens^F$, i.e.
\[ (a\la^F v)\ra^F b= (\partial F)^-(a_{[-1]}\tens v_{[-1]}\tens b_{[-1]}) \ 
a_{[0]}\la^F(v_{[0]}\ra^F b_{[0]}).\]
This all works in the same way as for the product $\mu^F$.  

\medskip
For $A$-bimodules $V$ and $W$, we can quotient the tensor product $V\tens W$ to get the 
tensor product over the algebra, $V\tens_A W$. The vector space which is quotiented by
is spanned by everything of the form $v\ra a\tens w-v\tens a\la w$. 
When our bimodules are in $({}^H\mathcal{M},\tens^F)$ we can proceed analogously by quotienting the vector space $V\tens^F W$ to give $V\tens_{A^F} W$
by identifying the elements given by going the two
ways round the following diagram:
\begin{eqnarray}\label{vcusakvcv}
\xymatrix{
V\tens^F(A\tens^F W) \ar[rr]_{\id\tens^F\la^F} \ar[dr]_{\Phi^{-1}} & & V\tens^F W \\
& (V\tens^F A)\tens^F W \ar[ur]_{\ra^F\tens^F\id} &
}
\end{eqnarray}
In fact this gives the same vector space $V\tens_A W$ as in the untwisted case. 
 The easiest way to see this is that the natural transformation $c$ between
the two tensor product functors, $\tens$ and $\tens^F$  (see \ref{kugdscvahksuv}) 
sends elements of the form
 $v.a\tens w-v\tens a.w$ exactly to the required relations for $V\tens_{A_F} W$, given by 
 (\ref{vcusakvcv}). Thus performing the quotient of vector spaces gives an isomorphism in the category. 

\medskip Now $V\tens_A W$ also has a standard $A$-bimodule structure, 
given by
\begin{eqnarray*}
a\la(v\tens w) \ =\ (a\la v)\tens w\ ,\quad (v\tens w)\ra a \ =\ v\tens (w\ra a)\ .
\end{eqnarray*}
Corresponding to this we have 
 left and right actions of $(A,\mu^F)$ on $V\tens^F W$ which are given by the following compositions:
\begin{eqnarray*}
&&A\tens^F (V\tens^F W) \stackrel{\Phi^{-1}} \longrightarrow (A\tens^F V)\tens^F W 
 \stackrel{\la^F\tens\id} \longrightarrow V\tens^F W\ ,\cr
 &&(V\tens^F W)\tens A \stackrel{\Phi} \longrightarrow V\tens^F (W\tens^F A) 
 \stackrel{\id\tens\ra^F} \longrightarrow V\tens^F W\ .
\end{eqnarray*}
The reader can check that these quotient to well defined maps on $V\tens_{A_F} W$.  These constructions make sense for algebras and their bimodules in any $k$-linear monoidal category, not just $({}^H\CM,\tens^F)$ and work the same way as for usual algebra due to Mac Lane's coherence theorem -- we need only express usual linear algebra constructions as compositions of maps and use the same in the monoidal category with the relevant $\Phi$ inserted to make sense.

\medskip
It will be useful to note that the definition of the left action on the tensor product (a similar picture holds for right)
is made so that the following diagram commutes:
\begin{eqnarray}\label{vcusfghgakvcv}
\xymatrix{
V \tens W    \ar[d]_{c} &   A \tens V \tens W \ar[r]_{c_{A,V\tens W}}   \ar[l]_{\la\tens\id}
& A \tens^F (V \tens W)  \ar[d]_{\id\tens^F c}   \\
V \tens^F W & &A \tens^F (V \tens^F W) \ar[ll]_{\la^F}
}
\end{eqnarray}
This just says that $\la^F$ is the image of $\la$ under the functorial equivalence from the category with $\tens$ to the category with $\tens^F$. When all maps are `quantized' by this equivalence then any equation expressed as a commutative diagram between maps will map over to an analogous commutative diagram on the $\tens^F$ side. This is similarly the abstract origin of the results in this paper.

\medskip
Such cochain twist ideas were applied to the octonions with $H=A=k(\Z_2^3)\isom k\Z_2^3$ and $\lambda=\Delta$ the `left regular coaction' of any group on itself. We use the second `group algebra' form where the algebra is spanned by basis elements $e_a$ where $a\in \Z_2^3$. The cochain $F$ and the twisted product of $A_F$ are
\begin{equation}\label{Focto}F(e_a\tens e_b)=(-1)^{\sum_{i\le j}a_ib_j+ a_1a_2b_3+a_1 b_2 a_3+b_1 a_2 a_3},\quad e_a\bullet e_b=F(e_a\tens e_b)e_{a+b}.\end{equation}
The product $\circ$ of $H$ is not, at it happens, changed (it remains associative) but this is largely an accident and its tensor category of comodules $\tens^F$ is nevertheless modified. Technically, it becomes a cotriangular coquasiHopf algebra which happens to be a Hopf algebra. The category with $\tens^F$ is a symmetric monodal one and in this category the octonions are both associative and commutative. Similarly if $C(G)$ is the algebra of functions for a finite or algebraic group and $F$ any choice of cochain on it, we have $C(G)^F$ with product $\circ$ which can happen to remain associative (e.g. in the case of Drinfelds $F$). We have a canonical coaction of $C(G)$ on itself by $\lambda=\Delta$ and hence a typically nonassociative algebra $C(G)_F$ with the induced  $\bullet$ product. Note the exact parallel with the octonion construction. They are both examples of systematic quantisation by cochain twist.

\section{Twisting of $\star$-structures}

We now proceed to study the different layers of geometry, but allowing our `coordinate algebras' $A$ to be potentially nonassociative by working in a monoidal category. The first item pertains to the right concept of real form or $*$-structure. In noncommutative geometry and in physics this specifies which elements are to represented as Hermitian or unitary in the sense that $*$ are to be respected. A categorical approch to this was introduced in \cite{BM3} with the notion of a `star object' in a bar category. We are only going to discuss it in the case relevant to cochain twists and we refer the reader to the paper for the formal category theory. 

\medskip
Thus we take $H$ to be a $*$-Hopf algebra, i.e. equipped with a conjugate-linear involution  
 $h\mapsto h^*$ on the algebra and required also to be a coalgebra map
 and order-reversing for the product. Then for every object $V$ in the comodule category ${}^H\mathcal{M}$  we have another $\overline  V$.
As a set $\overline{V}$ is the same as $V$, but we distinguish between elements 
by writing $\overline{v}\in \overline{V}$ for $v\in V$. Then the left $H$-coaction is given by
\begin{eqnarray*}
\lambda(\overline{v}) \,=\, (v_{[-1]})^* \tens \overline{v_{[0]}}\ .
\end{eqnarray*}
We have natural isomorphisms 
$\Upsilon:\overline{V\tens W}\to \overline{W}\tens\overline{V}$
and $\mathrm{bb}: V\to \overline{\overline{V}}$ defined by
\begin{eqnarray*}
\Upsilon(\overline{v\tens w}) \,=\, \overline{w} \tens \overline{v}\ ,\quad
\mathrm{bb}(v)\,=\, \overline{\overline{v}}\ .
\end{eqnarray*}
This provides the bar category associated to a Hopf $*$-algebra as a somewhat formal way of speaking about conjugate corepresentations.

\medskip
We now twist this construction by altering the map $\Upsilon$ to $\Upsilon^F
:\overline{V\tens^F W}\to \overline{W}\tens^F\overline{V}$ as follows:
\begin{eqnarray*}
\Upsilon^F(\overline{v\tens^F w}) \,=\, F(v_{[-2]}\tens w_{[-2]})^*\ 
F^-({w_{[-1]}}^*\tens {v_{[-1]}}^*)\ 
\overline{w_{[0]}} \tens^F \overline{v_{[0]}}\ ,\cr
\Upsilon^{F-1}(\overline{v}\tens^F \overline{w}) \,=\, F({v_{[-2]}}^*\tens {w_{[-2]}}^*)\ 
F^-({w_{[-1]}}\tens {v_{[-1]}})^*\ 
\overline{w_{[0]} \tens^F v_{[0]}}\ .
\end{eqnarray*}
These formulae are made to make the following diagram commute:
\begin{eqnarray*}
\xymatrix{
\overline{V\tens W}  \ar[r]^{\Upsilon} \ar[d]_{\overline{c}} 
& \overline{W}\tens\overline{V} \ar[d]_{c} \\
\overline{V\tens^F W}  \ar[r]^{\Upsilon^F}&
 \overline{W}\tens^F\overline{V}    }
\end{eqnarray*}

\medskip
\noindent To check that this works, we calculate
\begin{eqnarray*}
\Upsilon^F\,\overline{\Upsilon^F}\overline{\overline{v\tens^F w}} &=& F(v_{[-2]}\tens w_{[-2]})\ 
F^-({w_{[-1]}}^*\tens {v_{[-1]}}^*)^*\ 
\Upsilon^F(\overline{\overline{w_{[0]}} \tens^F \overline{v_{[0]}}}) \cr
&=& F(v_{[-4]}\tens w_{[-4]})\ 
F^-({w_{[-3]}}^*\tens {v_{[-3]}}^*)^*\ F({w_{[-2]}}^*\tens {v_{[-2]}}^*)^* \cr
&& F^-(v_{[-1]}\tens w_{[-1]})\, \overline{\overline{v_{[0]}}} \tens^F \overline{\overline{w_{[0]}}}\cr
&=& \overline{\overline{v}} \tens^F \overline{\overline{w}}\ .
\end{eqnarray*}
which is one of the axioms of a bar category. In this way we arrive at $({}^H\mathcal{M},\tens^F)$  as another bar category\cite{BM3}. It can be viewed as the category of comodules of $H^F$ as a certain $*$-coquasi-Hopf algebra. 

Next as in \cite{BM3}  we think of an $H$-comodule equipped with an antilinear involution categorically as a {\em star object}, meaning $V$ in ${}^H\mathcal{M}$ equipped with a morphism $\star_V:V\to \overline {V}$
so that $\overline{\star_V}\,\star_V=\mathrm{bb}_V$. We shall write $\star\, v=\overline{v^*}$ to relate it to the usual $*$ as an antilinear involution on a vector space $V$.  As $\star_V:V\to \overline {V}$ is a morphism
(i.e.\ a right $H$-comodule map) we deduce that
\begin{eqnarray*}
{v^*}_{[-1]}\tens \overline{{v^*}_{[0]}} \,=\, {v}_{[-1]}\tens \overline{{v_{[0]}}^*}\ .
\end{eqnarray*}
We can then define a star algebra $A$ in $({}^H\mathcal{M},\tens)$ as a star object $A$ which has an
associative product $\mu(a\tens b)=ab$ so that $\overline{\mu}\, \Upsilon^{-1}(\star\tens\star)=\star\mu$, which is a rather formal way to say $(ab)^*=b^*a^*$.  Putting our notions this way of course now makes sense in any bar category.

\begin{propos}
If $A$ is a star algebra $(A,\mu)$ in $({}^H\mathcal{M},\tens)$ then its cochain twist $(A,\mu^F)$ obeys 
$\overline{\mu^F}\, \Upsilon^{F-1}(\star\tens^F\star)=\star\mu^F$ and so is a star algebra in $({}^H\mathcal{M},\tens^F)$. Here the $\star$ operation itself is not deformed.
\end{propos}
\noindent {\bf Proof:}\quad By the definitions,
\begin{eqnarray*}
\Upsilon^{F-1}(\star\tens^F\star)(a\tens^F b) &=& \Upsilon^{F-1}(\overline{a^*}\tens^F \overline{b^*}) \cr
&=& F({a_{[-2]}}\tens {b_{[-2]}})\ 
F^-({b_{[-1]}}^*\tens {a_{[-1]}}^*)^*\ 
\overline{{b_{[0]}}^* \tens^F {a_{[0]}}^*}\ ,\cr
\overline{\mu^F}\, \Upsilon^{F-1}(\star\tens^F\star)(a\tens^F b) &=& 
F({a_{[-1]}}\tens {b_{[-1]}})\ \overline{{b_{[0]}}^*\, {a_{[0]}}^*} \cr
&=& F({a_{[-1]}}\tens {b_{[-1]}})\ \overline{(a_{[0]}b_{[0]})^*} \cr
&=& \star\mu^F(a\tens^F b)\ .\quad\square
\end{eqnarray*}

Implicit in this construction is the relation between $\star$ and the tensor product. 
That is, given that $V$ and $W$ are star objects, is $V\tens W$ a star object?
To simplify this we assume that
$H$ is commutative, and then we define $(v\tens w)^*=v^*\tens w^*$. 
We define $\Gamma:V\tens W\to W\tens V$ to be
$(\star_W^{-1}\tens \star_V^{-1})\,\Upsilon\,\star_{V\tens W}$ which merely recovers
usual transposition in the case of $({}^H\CM,\tens)$. In the twisted case, we define
\begin{eqnarray} \label{bvcadhlscv}
\star(v\tens^F w) \,=\, F(v_{[-2]}\tens w_{[-2]})\ F^-({v_{[-1]}}^*\tens {w_{[-1]}}^*)^*\ 
\overline{{v_{[0]}}^*\tens^F {w_{[0]}}^*}\ ,
\end{eqnarray}
and this corresponds to
\begin{eqnarray*}
\Gamma^F(v\tens^F w) &=& F(v_{[-2]}\tens w_{[-2]})\ F^-(w_{[-1]}\tens v_{[-1]})\ 
w_{[0]}\tens^F v_{[0]}\ .
\end{eqnarray*}
This is a generalised braiding associated to star objects and in some cases it obeys the braid relations \cite{BM3}. 

\medskip The reader might want an idea of real elements that is preserved under tensor product, that is that the tensor product of two real elements is real. So if we have $v=v^*$  and $w=w^*$,
then we would like $(v\tens^F w)^*=v\tens^F w$. From  (\ref{bvcadhlscv}) we can see how to ensure that this happens, we can insist on the following extra condition on the cochain:
\begin{eqnarray}
F(a\tens b)^*\ =\ F(a^*\tens b^*)\ .
\end{eqnarray}

\section{Twisting differential calculi and connections}
Suppose that $(\Omega A,\extd,\wedge)$ is a differential graded algebra in ${}^H\mathcal{M}$.
Here the exterior algebra is a direct sum of different degrees of differential form, with $\Omega^0 A=A$ itself, and there is a wedge product $\wedge$ and an exterior derivative $\extd:\Omega^n A\to \Omega^{n+1}A$ obeying $\extd^2=0$ and a graded Leibniz rule with respect to it.  We suppose that $H$ coacts on $A$ and that this coaction extends to the differential forms with $\extd$ and $\wedge$ respecting the coaction, or put another way, that all structure maps live in ${}^H\CM$. 

The differential calculus is twisted to $(\Omega A,\extd,\wedge^F)$, where
$\extd$ is unchanged, and
\begin{eqnarray*}
\xi\wedge^F\eta\,=\, F(\xi_{[-1]}\tens \eta_{[-1]})\ \xi_{[0]}\wedge \eta_{[0]}\ .
\end{eqnarray*}
We do not have to do any new work here, we just apply the cochain twist machinery of Section~2 to the algebra $\Omega A$ with its wedge product\cite{AkrMa}cf.\cite{MaOec}. 

Once we have a differential calculus we can define the notion of a connection or covariant derivative. For any usual associative algebra $A$ this is given by a map $\nabla:V\to \Omega^1 A\tens_A V$ where $V$ is a left $A$-module (more properly projective module for a vector bundle in noncommutative geometry), and required to obey
\begin{equation}\label{nabla}\nabla(av)=\extd a\tens_A v+a\nabla(v).\end{equation}
Its curvature is defined as a map $R_\nabla:V\to\Omega^2 A\tens_A V$ defined by 
\begin{equation}\label{Rnabla}R_\nabla=(\extd\tens\id-\id\wedge \nabla)\nabla=\nabla^{[1]}\nabla. \end{equation}
where the extension to higher forms is provided by 
\begin{eqnarray}\label{nablan} \nabla^{[n]}=\extd\tens\id + (-1)^n\id\wedge\nabla: \Omega^nA\tens _A V\to \Omega^{n+1}\tens_A V.\end{eqnarray}

\begin{propos} \label{vcdsuk}
Given a covariant derivative $\nabla$ on a left $A$-module $V$ which is also a left $H$-comodule map,
\begin{eqnarray*}
\nabla^F(v)\,=\, F^{-}(\xi_{[-1]}\tens w_{[-1]})\  \xi_{[0]}\tens^F w_{[0]}\end{eqnarray*}
where $\nabla(v)=\xi\tens w$ say (or a sum of such terms) provides a well-defined covariant derivative $\nabla^F:V\to\Omega^1A_F\tens_{A_F}V$ with $V$ twisted to an $A_F$-module as in Section~2.
\end{propos}
\noindent {\bf Proof:} We can say loosely that $\nabla^F=c\circ\nabla$ or that $\nabla^F$ is defined by commutativity of 
\begin{eqnarray}\label{vcakvcv}
\xymatrix{
V    \ar[r]_{\nabla^F}   \ar[d]^{\nabla}  &  \Omega^1 A \tens^F V  \\
 \Omega^1 A \tens V \ar[ur]_{c} & }
\end{eqnarray}
but we mean more precisely is choice of representatives for the result of $\nabla$ allowing us to work in $\Omega^1A\tens V$ followed by a check that $\nabla^F$ then descends to $\tens_{A_F}$ as required if $\nabla$ descends to $\tens_A$. A brief check will further show that this is a morphism in ${}^H\mathcal{M}$. To check the left
Liebniz rule, we begin by writing, for $v\in V$ and $b\in A$,
\begin{eqnarray*}
\nabla v &=& \xi\tens w\ ,\cr
v_{[-1]}\tens\nabla v_{[0]} &=& \xi_{[-1]}w_{[-1]}\tens \xi_{[0]}\tens w_{[0]}\ ,\cr
v_{[-1]}\tens\nabla(b. v_{[0]}) &=& \xi_{[-1]}w_{[-1]}\tens b.\xi_{[0]}\tens w_{[0]}
+v_{[-1]}\tens\extd b \tens v_{[0]}\ , \cr
v_{[-1]}\tens\nabla^F(b. v_{[0]}) &=& F^{-}(b_{[-1]}\xi_{[-1]}\tens w_{[-1]})\ \xi_{[-2]}w_{[-2]}\tens 
(b_{[0]}.\xi_{[0]}\tens^F w_{[0]}) \cr
&&
+\, F^{-}(b_{[-1]}\tens v_{[-1]})\ v_{[-2]}\tens(\extd b_{[0]} \tens^F v_{[0]})\ .
\end{eqnarray*}
Now we calculate, using the last equation, and writing $\bullet=\la^F$ for clarity,
\begin{eqnarray*}
\nabla^F(a\bullet v) &=& F(a_{[-1]}\tens v_{[-1]})\, \nabla^F(a_{[0]} v_{[0]}) \cr
&=& F(a_{[-2]}\tens \xi_{[-2]}w_{[-2]})\,F^{-}(a_{[-1]}\xi_{[-1]}\tens w_{[-1]})\ 
(a_{[0]}.\xi_{[0]}\tens^F w_{[0]}) \cr
&& +\, F(a_{[-2]}\tens v_{[-2]})\, F^{-}(a_{[-1]}\tens v_{[-1]})\ (\extd a_{[0]} \tens^F v_{[0]}) \cr
&=& F(a_{[-2]}\tens \xi_{[-2]}w_{[-2]})\,F^{-}(a_{[-1]}\xi_{[-1]}\tens w_{[-1]})\ 
(a_{[0]}.\xi_{[0]}\tens^F w_{[0]}) + \extd a\tens^F v\ .
\end{eqnarray*}
Next we calculate
\begin{eqnarray*}
a\bullet \nabla^F v &=& F^-(\xi_{[-1]}\tens w_{[-1]})\ a\bullet( \xi_{[0]}\tens^F w_{[0]}) \cr
&=& F^-(\xi_{[-2]}\tens w_{[-2]})\ (\partial F)(a_{[-1]}  \tens \xi_{[-1]}\tens w_{[-1]})
\ (a_{[0]} \bullet \xi_{[0]})\tens^F w_{[0]} \cr
&=& F^-(\xi_{[-3]}\tens w_{[-2]})\ (\partial F)(a_{[-2]}  \tens \xi_{[-2]}\tens w_{[-1]})\ 
F(a_{[-1]}  \tens \xi_{[-1]})
\ (a_{[0]} . \xi_{[0]})\tens^F w_{[0]} \cr
&=& F(a_{[-2]}\tens \xi_{[-2]}w_{[-2]})\,F^{-}(a_{[-1]}\xi_{[-1]}\tens w_{[-1]})\ 
(a_{[0]}.\xi_{[0]}\tens^F w_{[0]})\ .
\end{eqnarray*}
\quad $\square$

\medskip We next check what happens to maps which preserve the covariant derivatives
 between left modules. That is, given left modules and covariant derivatives
 $(V,\nabla_V)$ and  $(U,\nabla_U)$, we have a left module map
 $\theta:V\to U$ for which $(\id\tens\theta)\nabla_V=\nabla_U\,\theta:V\to \Omega^1 A\tens_A U$. 
 If we also suppose that $\theta:V\to U$ is an $H$-comodule map, then it is immediate from
 \ref{vcdsuk} that $(\id\tens^F\theta)\nabla^F_V=\nabla^F_U\,\theta:V\to \Omega^1 A\tens_{A_F} U$.
 We likewise extend the covariant derivative to maps
 \begin{eqnarray*}
\nabla^{F[n]}=\extd\tens^F\id+(-1)^n\id\wedge\nabla^F:\Omega^n A_F\tens_{A_F} V\to
 \Omega^{n+1} A_F\tens_{A_F} V\ ,
\end{eqnarray*}
where we remember that we have to use the associator to calculate $\id\wedge\nabla^F$. Then the curvature is defined analogously to before, as
\begin{eqnarray*}
R_{\nabla^F}=\nabla^{F[1]}\,\nabla^{F}:V\to \Omega^2 A_F\tens_{A_F} V\ .
\end{eqnarray*}

\begin{propos}\label{bvcshjkvk}
If the curvature of the covariant derivative $(V,\nabla)$ is given 
by $R_\nabla(v)=\omega\tens w\in \Omega^2 A\tens_A V$ (or a sum of such terms), then the 
curvature of the corresponding connection $(V,\nabla^F)$ is given by
\begin{eqnarray*}
R_{\nabla^F}(v)\,=\, F^-(\omega_{[-1]}\tens w_{[-1]})\  \omega_{[0]}\tens^F w_{[0]}\ .
\end{eqnarray*}
\end{propos}
\noindent {\bf Proof:}\quad 
We compute the curvature, taking $\nabla v = \xi\tens w$, as
\begin{eqnarray*}
\nabla^{F[1]}\nabla^F v &=& F^-(\xi_{[-1]}\tens w_{[-1]})\ \nabla^{F[1]}(\xi_{[0]}\tens^F w_{[0]}) \cr
&=& F^-(\xi_{[-1]}\tens w_{[-1]})\ \extd \xi_{[0]}\tens^F w_{[0]}
- F^-(\xi_{[-1]}\tens w_{[-1]})\ \xi_{[0]}\wedge^F \nabla^{F}w_{[0]}\ .
\end{eqnarray*}
Now we set $\nabla w = \eta\tens u$, and 
\begin{eqnarray*}
w_{[-1]}\tens\nabla^F( w_{[0]}) &=& F^{-}(\eta_{[-1]}\tens u_{[-1]})\ \eta_{[-2]}u_{[-2]}\tens 
(\eta_{[0]}\tens^F u_{[0]}) \cr
\end{eqnarray*}
Now we calculate
\begin{eqnarray*}
\nabla^{F[1]}\nabla^F v &=&
F^-(\xi_{[-1]}\tens w_{[-1]})\ \extd \xi_{[0]}\tens^F w_{[0]} \cr &&
-\ F^-(\xi_{[-1]}\tens \eta_{[-2]}u_{[-2]})\  
F^{-}(\eta_{[-1]}\tens u_{[-1]})\ \ 
\xi_{[0]}\wedge^F (\eta_{[0]}\tens^F u_{[0]}) \cr
&=&
F^-(\xi_{[-1]}\tens w_{[-1]})\ \extd \xi_{[0]}\tens^F w_{[0]} \cr &&
-\ F^-(\xi_{[-2]}\tens \eta_{[-3]}u_{[-3]})\  
F^{-}(\eta_{[-2]}\tens u_{[-2]})\ (\partial F)(\xi_{[-1]}\tens \eta_{[-1]}\tens u_{[-1]})
(\xi_{[0]}\wedge^F \eta_{[0]})\tens^F u_{[0]} \cr
&=&
F^-(\xi_{[-1]}\tens w_{[-1]})\ \extd \xi_{[0]}\tens^F w_{[0]} 
- F^-(\xi_{[-1]} \eta_{[-1]}\tens u_{[-1]})\  
(\xi_{[0]}\wedge \eta_{[0]})\tens^F u_{[0]} \ .\end{eqnarray*}
We compare this with the straight computation of $R_\nabla(v)$ and change notations to the form in the statement. We can say loosely that $R_{\nabla^F}=c\circ R_\nabla$. \quad $\square$

\section{Bimodule covariant derivatives}

In applying these ideas to Riemannan geometry we will be in the case of $V=\Omega^1 A$ which is a bimodule.  When $V$ is a bimodule we may have a well-defined  generalised braiding $\sigma: V\tens_A \Omega^1 A\to \Omega^1 A\tens_A V$ such that
\[ \nabla(v\ra a)=(\nabla v)\ra a+\sigma(v\tens_A\extd a).\]
Such an object when it exists is called a {\em bimodule covariant derivative}.  We already saw in Section~2 how to twist an $H$-covariant bimodule to one of $(A,\mu^F)$. We suppose that $\nabla$ and $\sigma$ are also $H$-covariant i.e. expressed in the category ${}^H\CM$.

\begin{propos}\label{twistnabla}
If $(V,\nabla,\sigma)$ is an $H$-covariant bimodule covariant derivative then 
$(V,\nabla^F,\sigma^F)$ is a bimodule covariant derivative for $(A,\mu^F)$ where
(writing $\sigma(v_{[0]}\tens\xi_{[0]})=\eta\tens w$ or a sum of such representatives)
\begin{eqnarray*}
\sigma^F(v\tens^F\xi) &=& F(v_{[-1]}\tens \xi_{[-1]})\  F^-\big(
\eta_{[-1]}\tens w_{[-1]} \big) \eta_{[0]}\tens^F w_{[0]} \ 
\end{eqnarray*}
descends to the quotients as required.
\end{propos}
\noindent {\bf Proof:}\quad One should either work with $\tens_A$ and $\tens_{A_F}$ or representatives with $\tens$ and $\tens^F$ and an eventual quotient understood. By definition,
\begin{eqnarray*}
\sigma^F(v\tens^F \extd a) &=& \nabla^F(v\ra^F a)-\nabla^F(v)\ra^F a\cr
&=& F(v_{[-1]}\tens a_{[-1]})\ \nabla^F(v_{[0]}\,a_{[0]})-\nabla^F(v)\ra^F a\ .
\end{eqnarray*}
If we write $\nabla v=\xi\tens w$, then $v_{[-1]}\tens\nabla v_{[0]}=
\xi_{[-1]}\, w_{[-1]}\tens \xi_{[0]}\tens w_{[0]}$, so
\begin{eqnarray*}
\sigma^F(v\tens^F \extd a) &=& 
F(v_{[-1]}\tens a_{[-1]})\ \nabla^F(v_{[0]}\,a_{[0]})-
F^-(\xi_{[-1]}\tens w_{[-1]})\ (\xi_{[0]}\tens^F w_{[0]})\ra^F a\ .
\end{eqnarray*}
Now we use
\begin{eqnarray*}
v_{[-1]}\tens a_{[-1]} \tens\nabla(v_{[0]} a_{[0]} ) &=& v_{[-1]}\tens a_{[-1]} \tens\nabla(v_{[0]}) a_{[0]}  +
v_{[-1]}\tens a_{[-1]}  \tens\sigma(v_{[0]}\tens\extd a_{[0]} ) \cr
&=& 
\xi_{[-1]}\, w_{[-1]}\tens a_{[-1]} \tens \xi_{[0]}\tens w_{[0]} a_{[0]}  +
v_{[-1]}\tens a_{[-1]}  \tens\sigma(v_{[0]}\tens\extd a_{[0]} ) \ ,\cr
v_{[-1]}\tens a_{[-1]} \tens\nabla^F(v_{[0]} a_{[0]} ) &=& 
\xi_{[-2]}\, w_{[-2]}\tens a_{[-2]} \tens F^-( \xi_{[-1]}\tens w_{[-1]} a_{[-1]})\ ( \xi_{[0]}\tens^F w_{[0]} a_{[0]})
\cr &&+\,
v_{[-1]}\tens a_{[-1]}  \tens 
F^-(\sigma^{(1)}(v_{[0]}\tens\extd a_{[0]} )_{[-1]}\tens \sigma^{(2)}(v_{[0]}\tens\extd a_{[0]} )_{[-1]})\cr
&&
\sigma^{(1)}(v_{[0]}\tens\extd a_{[0]} )_{[0]}\tens \sigma^{(2)}(v_{[0]}\tens\extd a_{[0]} )_{[0]}\ .
\end{eqnarray*}
The first term in this equation cancels with $\nabla^F(v)\ra^F a$, giving the answer.
\quad$\square$

\medskip The main purpose of the bimodule covariant derivative construction 
is to allow us to tensor product covariant derivatives. Given $A$-bimodule covariant derivatives
$(V,\nabla_V,\sigma_V)$ and $(W,\nabla_W,\sigma_W)$, we have the following covariant derivative on 
$V\tens_A W$ defined by
\begin{eqnarray}\label{vcusakvcvxx}
\nabla_{V\tens_A W}(v\tens w) &=& \nabla_V(v)\tens w + (\sigma_V\tens_A\id)(v\tens\nabla_W w)\ .
\end{eqnarray}
quotiented to $V\tens_AW$. By Proposition~\ref{twistnabla} we have a bimodule connection $\nabla^F$ on $V\tens_AW$ as an $A_F$-bimodule $V\tens_{A_F}W$. 

\begin{propos}
Given $A$-bimodule covariant derivatives
$(V,\nabla_V,\sigma_V)$ and $(W,\nabla_W,\sigma_W)$,
the twisted $\nabla^F$ on $V \tens_{A_F} W$ obeys
\begin{eqnarray*}
\nabla_{V\tens_{A_F} W}^F \ =\ \Phi(\nabla_{V}^F\tens^F \id) +(\sigma_V^F\tens^F \id)\Phi^{-1}
(\id\tens^F \nabla_{W}^F))\ .
\end{eqnarray*}
where the composition on the right descends to the required quotient over $A_F$. \end{propos}
\noindent {\bf Proof:}\quad Unwinding the definitions here, $\nabla^F$ is characterised by\begin{eqnarray}\label{vcusakvcvyy*}
\xymatrix{
V\tens W \ar[rr]^{c} \ar[d]_{\nabla} & & V\tens^F W  \ar[d]_{\nabla^F} \\
\Omega^1 A \tens V \tens W \ar[r]^{c}&
\Omega^1 A \tens^F (V \tens W) \ar[r]^{\id \tens^F c} &
\Omega^1 A \tens^F (V \tens^F W)}
\end{eqnarray}
where we again abuse notations as in previous proofs. To check that the map is indeed well defined over $V \tens_{A_F} W$,
we check that the paths from the top left to bottom right corners of the following diagrams are identical (up to the relations on $V \tens_{A_F} W$ in the bottom right corner), which can be done easily just by considering the left and bottom paths of the diagrams:
\begin{eqnarray*}
\xymatrix{
V\tens A\tens W  \ar[d]_{\ra\tens\id} \ar[r]_{c}  & (V\tens A)\tens^F W  \ar[r]_{c\tens^F\id}  & 
(V\tens^F A)\tens^F W  \ar[d]_{\ra^F\tens^F\id} \\
V\tens W \ar[rr]^{c} \ar[d]_{\nabla} & & V\tens^F W  \ar[d]_{\nabla^F} \\
\Omega^1 A \tens V \tens W \ar[r]^{c}&
\Omega^1 A \tens^F (V \tens W) \ar[r]^{\id \tens^F c} &
\Omega^1 A \tens^F (V \tens^F W)}
\end{eqnarray*}
\begin{eqnarray*}
\xymatrix{
V\tens A\tens W  \ar[d]_{\id\tens\la} \ar[r]_{c}  & V\tens^F (A\tens W)  \ar[r]_{\id\tens^F c}  & 
V\tens^F (A\tens^F W)  \ar[d]_{\id\tens^F\la^F} \\
V\tens W \ar[rr]^{c} \ar[d]_{\nabla} & & V\tens^F W  \ar[d]_{\nabla^F} \\
\Omega^1 A \tens V \tens W \ar[r]^{c}&
\Omega^1 A \tens^F (V \tens W) \ar[r]^{\id \tens^F c} &
\Omega^1 A \tens^F (V \tens^F W)}
\end{eqnarray*}
We do not actually need to check this nor that we have a bimodule connection, in view of Proposition~\ref{twistnabla}. However, it is instructive to do so for good measure. For the left Liebnitz rule, we use
\begin{eqnarray}\label{cvadsjkcv}
\xymatrix{
A\tens V\tens W  \ar[d]_{\la\tens\id} \ar[r]_{c}  & A\tens^F (V\tens W)  \ar[r]_{\id\tens^F c}  & 
A\tens^F (V\tens^F W)  \ar[d]_{\la^F} \\
V\tens W \ar[rr]^{c} \ar[d]_{\nabla} & & V\tens^F W  \ar[d]_{\nabla^F} \\
\Omega^1 A \tens V \tens W \ar[r]^{c}&
\Omega^1 A \tens^F (V \tens W) \ar[r]^{\id \tens^F c} &
\Omega^1 A \tens^F (V \tens^F W)}
\end{eqnarray}
Summarising the right side of diagram (\ref{cvadsjkcv}) as `right', we get two parts to 
(\ref{cvadsjkcv}), the first being
\begin{eqnarray*}
\xymatrix{
A\tens V\tens W  \ar[d]_{\extd\tens\id\tens\id} \ar[r]_{c}  & A\tens^F (V\tens W)  \ar[r]_{\id\tens^F c}  & 
A\tens^F (V\tens^F W)  \ar[d]_{\mathrm{first\ right}} \\
\Omega^1 A \tens V \tens W \ar[r]^{c}&
\Omega^1 A \tens^F (V \tens W) \ar[r]^{\id \tens^F c} &
\Omega^1 A \tens^F (V \tens^F W)}
\end{eqnarray*}
which gives the right hand side $\extd\tens^F\id$ by functoriality. The second term is
rather more difficult:
\begin{eqnarray*}
\xymatrix{
A\tens (V\tens W)  \ar[d]_{\id\tens\nabla} \ar[r]_{c}  & A\tens^F (V\tens W)  \ar[r]_{\id\tens^F c}  & 
A\tens^F (V\tens^F W)   \ar[dd]_{\mathrm{second\ right}} \\
A\tens \Omega^1 A \tens V\tens W  \ar[d]_{\la\tens\id\tens\id} & & \\
\Omega^1 A \tens V \tens W \ar[r]^{c}&
\Omega^1 A \tens^F (V \tens W) \ar[r]^{\id \tens^F c} &
\Omega^1 A \tens^F (V \tens^F W)}
\end{eqnarray*}
Functoriality of $c$ applied to the top and bottom sides gives
\begin{eqnarray*}
\xymatrix{
A\tens (V\tens W)  \ar[d]_{\id\tens\nabla} \ar[r]_{\id\tens c}  & A\tens (V\tens^F W)  \ar[r]_{c}  & 
A\tens^F (V\tens^F W)   \ar[dd]_{\mathrm{second\ right}} \\
A\tens \Omega^1 A \tens V\tens W  \ar[d]_{\la\tens\id\tens\id} & & \\
\Omega^1 A \tens V \tens W \ar[r]^{\id\tens c}&
\Omega^1 A \tens (V \tens^F W) \ar[r]^{c} &
\Omega^1 A \tens^F (V \tens^F W)}
\end{eqnarray*}
and since the bottom left two operators commute, we get
\begin{eqnarray*}
\xymatrix{
A\tens (V\tens W)  \ar[d]_{\id\tens\nabla} \ar[r]_{\id\tens c}  & A\tens (V\tens^F W)  \ar[r]_{c}  & 
A\tens^F (V\tens^F W)   \ar[dd]_{\mathrm{second\ right}} \\
A\tens \Omega^1 A \tens V\tens W  \ar[d]_{\id\tens\id\tens c} & & \\
A\tens \Omega^1 A \tens (V \tens^F W) \ar[r]^{\la\tens \id}&
\Omega^1 A \tens (V \tens^F W) \ar[r]^{c} &
\Omega^1 A \tens^F (V \tens^F W)}
\end{eqnarray*}
Now use diagram (\ref{vcusakvcvyy*}) to get
\begin{eqnarray*}
\xymatrix{
A\tens (\Omega^1 A \tens^F(V\tens^F W))  \ar[d]_{\id\tens c^{-1}}   &
 A\tens (V\tens^F W)  \ar[r]_{c}  \ar[l]_{\id\tens\nabla^F}& 
A\tens^F (V\tens^F W)   \ar[d]_{\mathrm{second\ right}} \\
A\tens \Omega^1 A \tens (V \tens^F W) \ar[r]^{\la\tens \id}&
\Omega^1 A \tens (V \tens^F W) \ar[r]^{c} &
\Omega^1 A \tens^F (V \tens^F W)}
\end{eqnarray*}
More functoriality of $c$ gives
\begin{eqnarray*}
\xymatrix{
A\tens^F (\Omega^1 A \tens(V\tens^F W)) \ar[d]_{c^{-1}}   &
A\tens^F (\Omega^1 A \tens^F(V\tens^F W))  \ar[l]_{\id\tens c^{-1}} & 
A\tens^F (V\tens^F W)  \ar[l]_{\id\tens^F\nabla^F}  \ar[d]_{\mathrm{second\ right}} \\
A\tens \Omega^1 A \tens (V \tens^F W) \ar[r]^{\la\tens \id}&
\Omega^1 A \tens (V \tens^F W) \ar[r]^{c} &
\Omega^1 A \tens^F (V \tens^F W)}
\end{eqnarray*}
Now a direct calculation will show that this is the required
$\la^F\,(\id\tens^F\nabla^F:A\tens^F (V\tens^F W)\to \Omega^1 A \tens^F (V \tens^F W)$.

For the bimodule part of the bimodule covariant derivative, consider
\begin{eqnarray*}
\xymatrix{
V\tens W\tens A  \ar[d]_{\id\tens\ra} \ar[r]_{c}  & (V\tens W)\tens^F A  \ar[r]_{c\tens^F\id}  & 
(V\tens^F W)\tens^F A  \ar[d]_{\ra^F} \\
V\tens W \ar[rr]^{c} \ar[d]_{\nabla} & & V\tens^F W  \ar[d]_{\nabla^F} \\
\Omega^1 A \tens V \tens W \ar[r]^{c}&
\Omega^1 A \tens^F (V \tens W) \ar[r]^{\id \tens^F c} &
\Omega^1 A \tens^F (V \tens^F W)}
\end{eqnarray*}
This splits into two terms, the first of which is
\begin{eqnarray*}
\xymatrix{
V\tens W\tens A  \ar[d]_{\id\tens\id\tens\extd} \ar[r]_{c}  & (V\tens W)\tens^F A  \ar[r]_{c\tens^F\id}  & 
(V\tens^F W)\tens^F A  \ar[dd]_{\mathrm{first\ term}} \\
V\tens W \tens\Omega^1 A \ar[d]_{\sigma_{V\tens W}} & & \\
\Omega^1 A \tens V \tens W \ar[r]^{c}&
\Omega^1 A \tens^F (V \tens W) \ar[r]^{\id \tens^F c} &
\Omega^1 A \tens^F (V \tens^F W)}
\end{eqnarray*}
By functoriality of $\extd$ and the formula for $\sigma_{V\tens W}$, this is
\begin{eqnarray*}
\xymatrix{
V\tens W\tens \Omega^1 A  \ar[d]_{\id\tens\sigma_W}   &
 (V\tens W)\tens^F \Omega^1 A \ar[l]_{c^{-1}} & 
(V\tens^F W)\tens^F A  \ar[dd]_{\mathrm{first\ term}}  \ar[l]_{c^{-1}\tens^F\extd}  \\
V\tens  \Omega^1 A \tens W  \ar[d]_{\sigma_{V}\tens \id} & & \\
\Omega^1 A \tens V \tens W \ar[r]^{c}&
\Omega^1 A \tens^F (V \tens W) \ar[r]^{\id \tens^F c} &
\Omega^1 A \tens^F (V \tens^F W)}
\end{eqnarray*}
The second term is
\begin{eqnarray*}
\xymatrix{
V\tens W\tens A  \ar[d]_{\nabla_{V\tens W}\tens\id} \ar[r]_{c}  & (V\tens W)\tens^F A  \ar[r]_{c\tens^F\id}  & 
(V\tens^F W)\tens^F A  \ar[dd]_{\mathrm{second\ term}} \\
\Omega^1 A \tens V\tens W \tens A\ar[d]_{\id\tens\id\tens\ra} & & \\
\Omega^1 A \tens V \tens W \ar[r]^{c}&
\Omega^1 A \tens^F (V \tens W) \ar[r]^{\id \tens^F c} &
\Omega^1 A \tens^F (V \tens^F W)}
\end{eqnarray*}
By using the functoriality of $c$, this is
\begin{eqnarray*}
\xymatrix{
(\Omega^1 A \tens V\tens W)\tens^F A  \ar[d]_{c^{-1}}  &
 (V\tens W)\tens^F A  \ar[r]_{c\tens^F\id} \ar[l]_{\nabla_{V\tens W}\tens^F\id}   & 
(V\tens^F W)\tens^F A  \ar[dd]_{\mathrm{second\ term}} \\
\Omega^1 A \tens V\tens W \tens A\ar[d]_{\id\tens\id\tens\ra} & & \\
\Omega^1 A \tens V \tens W \ar[r]^{c}&
\Omega^1 A \tens^F (V \tens W) \ar[r]^{\id \tens^F c} &
\Omega^1 A \tens^F (V \tens^F W)}
\end{eqnarray*}

The first term of the equation in the statement is given by differentialting the $v$, and can be read off the diagram as the following, where we set
$\nabla v=\xi\tens u$:
\begin{eqnarray*}
F(\xi_{[-2]} u_{[-3]} \tens w_{[-3]})\,F^-(\xi_{[-1]} u_{[-2]} \tens w_{[-2]})\,
F^-(u_{[-1]} \tens w_{[-1]})\,\xi_{[0]}\tens^F(u_{[0]}\tens^F w_{[0]})\ ,
\end{eqnarray*}
which gives $\Phi(\nabla_{V}^F\tens^F \id)$. To find the second term,
we differentiale the $w$ and substitute into diagram (\ref{vcusakvcvyy*}) to get
\begin{eqnarray*}\label{vcusakvcvyyyy}
\xymatrix{
V\tens W \ar[rr]^{c} \ar[d]_{\id\tens\nabla_W} & & V\tens^F W  \ar[d]_{\mathrm{second}} \\
V\tens \Omega^1 A \tens W  \ar[d]_{\sigma\tens\id} & 
 & \Omega^1 A \tens^F (V \tens^F W) \\
\Omega^1 A \tens V \tens W \ar[r]^{c}& 
(\Omega^1 A \tens V) \tens^F W \ar[r]^{c \tens^F \id} &
(\Omega^1 A \tens^F V) \tens^F W    \ar[u]^{\Phi}  }
\end{eqnarray*}
By the functoriality of $c$ we get
\begin{eqnarray*}\label{vcusakvcvyyyy1}
\xymatrix{
V\tens W \ar[rr]^{c} \ar[d]_{\id\tens\nabla_W} & & V\tens^F W  \ar[d]_{\mathrm{second}} \\
V\tens \Omega^1 A \tens W  \ar[d]_{c} & 
 & \Omega^1 A \tens^F (V \tens^F W) \\
(V\tens \Omega^1 A)  \tens^F W \ar[r]^{\sigma\tens^F\id}& 
(\Omega^1 A \tens V) \tens^F W \ar[r]^{c \tens^F \id} &
(\Omega^1 A \tens^F V) \tens^F W    \ar[u]^{\Phi}  }
\end{eqnarray*}
and by the definition of $\sigma^F$,
\begin{eqnarray*}\label{vcusakvcvyyyy2}
\xymatrix{
V\tens W \ar[rr]^{c} \ar[d]_{\id\tens\nabla_W} & & V\tens^F W  \ar[d]_{\mathrm{second}} \\
V\tens \Omega^1 A \tens W  \ar[d]_{c} & 
 & \Omega^1 A \tens^F (V \tens^F W) \\
(V\tens \Omega^1 A)  \tens^F W \ar[r]^{c \tens^F \id}& 
(V\tens^F \Omega^1 A) \tens^F W \ar[r]^{\sigma^F\tens^F\id} &
(\Omega^1 A \tens^F V) \tens^F W    \ar[u]^{\Phi}  }
\end{eqnarray*}
Finally we use the definition of $\Phi^{-1}$ and $\nabla^F_W$.\quad$\square$

\medskip
Finally, we need to see how the notion of $\star$-compatibilty behaves. This is an important extra condition in noncommutative Riemannian geometry \cite{BM4} and corresponds in our application to ensuring reality constraints. To keep things simple  we again assume that the Hopf algebra $H$ is now commutative and has $S^2=\id$, and that
$(A,\Omega^*A,\extd)$ is a classical differential geometry. That means that
$A$ is commutative, for bimodule $A$-covariant derivatives
$\sigma$ is transposition, and that $\wedge$ is given by antisymmetrisation. 

\begin{propos}
Suppose for simplicity that we start in the classical situation with $A$ commutative. If $\nabla$ is compatible with $\star$, then $\nabla^F$ is compatible with $\star$. 
\end{propos}
\noindent {\bf Proof:}\quad In classical differential geometry, the condition that $\nabla$ is compatible with $\star$ is 
\begin{eqnarray*}
\overline{\sigma}\,\Upsilon^{-1}(\star\tens\star)\nabla(v) \,=\, \overline{\nabla(v^*)}\ ,
\end{eqnarray*}
and this reduces to $\nabla (v^*)=\xi^*\tens w^*$ where $\nabla (v)=\xi\tens w$. Now
\begin{eqnarray*}
\overline{\sigma^F}\,\Upsilon^{F-1}(\star\tens^F\star)\nabla^F v &=& 
F^-(\xi_{[-1]}\tens w_{[-1]})\, \overline{\sigma^F}\,\Upsilon^{F-1}(\star\tens^F\star)
(\xi_{[0]}\tens^F w_{[0]}) \cr
&=& F^-(\xi_{[-1]}\tens w_{[-1]})\, \overline{\sigma^F}\,\Upsilon^{F-1}
(\overline{{\xi_{[0]}}^*}\tens^F \overline{{w_{[0]}}^*}) \cr
&=& F^-(w_{[-1]}^* \tens \xi_{[-1]}^*)^* \, \overline{\sigma^F}\,
(\overline{{w_{[0]}}^*\tens^F {\xi_{[0]}}^*}) \cr
&=&F^-(w_{[-3]}^* \tens \xi_{[-3]}^*)^* \, F({w_{[-2]}}^* \tens {\xi_{[-2]}}^*)^* \, 
F^-({\xi_{[-1]}}^* \tens {w_{[-1]}}^*)^* (\overline{{\xi_{[0]}}^*\tens^F {w_{[0]}}^*}) \cr
&=& F^-({\xi_{[-1]}}^* \tens {w_{[-1]}}^*)^* (\overline{{\xi_{[0]}}^*\tens^F {w_{[0]}}^*}) \cr
&=& \overline{\nabla^F (v^*)}\ .\quad\square
\end{eqnarray*}

\section{Twisting of metrics and Riemannian structures}

We are now ready to put much of this together into an algebraic framework for Riemannian geometry. In this section we will use differential calculus over an (associative) algebra $A$, but the example to bear in mind is of the smooth functions on a manifold. As previously,
we take a differential calculus $\Omega A$ in the category ${}^H\CM$, where (because of our previous simplifications) we assume that $H$ is commutative. We again assume that the connection $\nabla$ to be deformed is a $H$-comodule map. To fit the twisting procedure, we assume that the Riemannian metric is $H$-invariant.
Within the standard framework of real differential geometry, star preservation is automatic
(as star is the identity), so we assume that
the connection preserves star.

\medskip
We look at connections $\nabla:\Omega^1 A\to \Omega^1 A\tens_A\Omega^1 A$ where $\Omega A=\oplus_n\, \Omega^n A$ is the exterior algebra on $A$.  We use the formalism invented for (associative) noncommutative geometry but we are mainly interested in using it in the classical commutative case and then  twisting it to generate nonassociative examples. Thus, the curvature $R_\nabla$ now has the meaning of Riemann curvature and $\nabla$ will be some kind of Levi-Civita connection for a metric. Torsion for any connection on $\Omega^1 A$ is defined as
\begin{eqnarray}\label{Tnabla}
{\rm Tor}_\nabla=\extd-\wedge\nabla: \Omega^1 A \to \Omega^2 A
\end{eqnarray}
so torsion-free makes sense.

To formulate the metric and metric compatibility there are currently two approaches. The original one\cite{Ma:rie}
is
\begin{eqnarray}\label{g}
g\in \Omega^1 A\tens_A\Omega^1 A,\quad \wedge(g)=0
\end{eqnarray}
where the second expresses symmetry with respect to the notion of skew-symmetrization used in the exterior algebra. We also require some form of non-degeneracy best expressed in terms of projective modules or in a frame bundle approach.  Finally, this setting has a remarkable symmetry between tangent and cotangent bundles and dual to torsion (i.e. torsion on the dual bundle regarded as cotangent bundle) is a notion of {\em cotorsion}. If we assume that the torsion vanishes then cotorsion reduces to
\begin{eqnarray}\label{coTnabla} 
\mathrm{coTor}\ =\ (\wedge\nabla\tens\id-(\wedge\tens\id)(\id\tens\nabla))g\in
\Omega^2 A\tens_A\Omega^1 A\ 
\end{eqnarray}
and we see that its vanishing is a weaker (skew-symmetrized) version of usual metric compatibility.  A generalised Levi-Civita connection is then a torsion free cotorsion free (or skew-compatible) connection. This weakening of the usual notion of metric compatibility seems to be necessary in several examples where the generalised Levi-Civita connection then exists and is unique for the chosen metric. If $A$ is commutative we can of course impose $\nabla g=0$ in the usual way extending it as a derivation to the two tensor factors but keeping the `left output' of $\nabla$ to the far left where it can be evaluated against a vector field.

A second most recent approach\cite{BM4} makes sense for bimodule connections. First of all, there is a useful weaker notion to torsion free, namely that the torsion ${\rm Tor}_\nabla$ be a bimodule map (we call this torsion-compatible).   We also make use of bar categories and work with Hermitian metrics $g\in \Omega^1 A\tens_A\overline{\Omega^1 A}$ and in this case full metric compatibility $\nabla g=0$ makes perfect sense. There is now a slightly different notion of cotorsion in place of (\ref{coTnabla}),
\begin{equation}\label{newcoT}
\mathrm{coTor}\ =\ (\wedge\tens\id)\nabla g\in
\Omega^2 A\tens_A \overline{\Omega^1 A}\ .
\end{equation}
 In \cite{BM4} we showed that these notions lead, for example, to a unique torsion free metric and $*$-compatible connection with classical limit on a class of metrics on $C_q(SU_2)$ with its left covariant 3D calculus.

When all of these structures are $H$-covariant, we can twist them. We have already seen how $\wedge,\nabla$ twist and checking our formulae for $\wedge^F$ and $\nabla^F$ we see that 
\[ {\rm Tor}_{\nabla^F}=\extd-\wedge^F\nabla^F={\rm Tor}_\nabla\]
 is unchanged as a linear map $\Omega^1 A \to \Omega^2 A$.  For $g\in \Omega^1 A\tens_A\Omega^1 A$, we take 
\[ g^F=c(g)\in \Omega^1 A_F\tens_{A_F}\Omega^1 A_F\]
and the cotorsion or skew-metric compatibility computed in the category with $\tens^F$ takes the form
\begin{eqnarray*}
\mathrm{coTor}^F\ =\ (\wedge^F\nabla^F\tens^F\id-(\wedge^F\tens^F\id)\Phi^{-1}(\id\tens^F\nabla^F))\,g^F \in
\Omega^2 A_F\tens_{A_F}\Omega^1 A_F\ 
\end{eqnarray*}
where an associator is inserted. The reader can easily check the following commutative diagram,
\begin{eqnarray*}
\xymatrix{
\Omega^1 A\tens\Omega^1 A \ar[r]^{c}  \ar[d]^{X} & \Omega^1 A\tens^F\Omega^1 A  \ar[d]^{X^F} \\
\Omega^2 A\tens\Omega^1 A \ar[r]^{c} & \Omega^2 A\tens^F\Omega^1 A
}
\end{eqnarray*}
where $X=(\wedge\tens\id)(\id\tens\nabla)$ and
$X^F=(\wedge^F\tens^F\id)\Phi^{-1}(\id\tens^F\nabla^F)$.
It follows that the twisted cotorsion vanishes if and only if the original cotorsion vanishes. Similarly in the hermitian metric framework. The following result follows directly from this discussion:

\begin{propos}  
Let 
$g\in \Omega^1 A\tens_A\Omega^1 A$ be a Riemannian metric $A$ invariant under $H$ and assume for simplicity that we are in the classical situation. Then a torsion free, cotorsion free connection $\nabla$ on $A$ which is an $H$-comodule map is twisted to a star preserving, torsion free, cotorsion free connection $\nabla^F$ on $A_F$.
\end{propos}

In particular, the classical Levi-Civita connection will be covariant and hence twist as required. We can also consider a Riemannian metric as an inner product
$\<,\>:\Omega^1 A\tens_A \overline{\Omega^1 A}\to A$
(using non-degeneracy), which we assume is an $H$-comodule map. We can use a little more generality to get the following result
for bimodule connections:

\begin{lemma}
Begin with $A$-bimodule connections $(V,\nabla_V)$ and $(W,\nabla_W)$
in ${}^H\mathcal{M}$.
Suppose that we have a map $\kappa:V\tens_A W\to A$ in ${}^H\mathcal{M}$ which preserves the connections (using $(A,\extd)$ as a connection on $A$), i.e.\
\begin{eqnarray*}
\extd\,\kappa\ =\ (\id\tens_A \kappa)\nabla_{V\tens_A W}:V\tens W\to \Omega^1 A\ .
\end{eqnarray*}
Then for the twisted connections we have
\begin{eqnarray*}
\extd\,\kappa^F\ =\ (\id\tens_{A_F}\kappa^F)\nabla^F_{V\tens_{A_F} W}:V\tens_{A_F} W\to \Omega^1 A_F\ ,
\end{eqnarray*}
where
\begin{eqnarray*}
 \kappa^F(v\tens^F w) \ =\ F(v_{[-1]}\tens w_{[-1]})\,\kappa(v_{[0]}\tens w_{[0]})\ 
\end{eqnarray*}
descends to $\kappa^F:V\tens_{A_F} W\to A_F$. 
\end{lemma}
\noindent {\bf Proof:}\quad From (\ref{vcusakvcvyy*}) we have
(where $\mu$ is multiplication)
\begin{eqnarray*}
\xymatrix{
&V\tens W \ar[rr]^{c} \ar[d]_{\nabla} \ar[dl]_{\kappa} & & V\tens^F W  \ar[d]_{\nabla^F}  \\
A \ar[d]_{\extd} &
\Omega^1 A \tens V \tens W \ar[r]^{c}\ar[d]^{\id\tens\kappa}&
\Omega^1 A \tens^F (V \tens W) \ar[r]^{\id \tens^F c} \ar[d]^{\id\tens^F\kappa}&
\Omega^1 A \tens^F (V \tens^F W) \ar[dl]^{\id\tens^F\kappa^F} \\
\Omega^1 A  &
\Omega^1 A \tens A  \ar[r]^{c}   \ar[l]^{\mu} & \Omega^1 A \tens^F A  &}
\end{eqnarray*}
The result can be read off from the diagram.\quad$\square$

\begin{cor}
Let  $g\in \Omega^1 A\tens_A\overline{\Omega^1 A}$ be a hermitian Riemannian metric $A$ invariant under $H$ and assume for simplicity that we are in the classical situation.   Then a torsion free bimodule connection $\nabla$ on $A$ which preserves the metric and is an $H$-comodule map is twisted to a star preserving,  torsion free bimodule connection $\nabla^F$ on $A_F$ which preserves the twisted Riemannian metric.
\end{cor}

Again, we can take a classical connection which is trivially a bimodule connection with the flip map for $\sigma$ and will obtain a twisted one.  Clearly one can go on in this line. One can use a metric to define interior product of 1-forms on $n$-forms and use this to define all relevant quantities in physics including candidates for the Ricci tensor and the stress energy tensor of scalar and other fields.  The problem at this point is not that we can't do this, but that there may be several different covariant formalisations all with the same classical limit and they will all twist. For this one really needs a theory that works beyond the twist examples. Until then, twisting nevertheless provides a source of examples.

This is illustrated by the Ricci tensor where the classical `formula' of contracting Riemann could be expressed in different ways depending on where we contract. While not unique, probably the simplest formulation is as follows, as a version similar to \cite{Ma:rie,Ma:sph}. First,  we need a splitting map such that
\begin{equation}\label{lift} i:\Omega^2 A\to \Omega^1 A\tens_A \Omega^1 A,\quad \wedge\circ i=\id. \end{equation}
In the classical case this is  $i(\xi\wedge\eta)=(\xi\tens\eta-\eta\tens\xi)/2$. Then we can define Ricci in our first formulation of metric as
\begin{equation}\label{ricci} {\rm Ricci}=(\<\ ,\ \>\tens\id)(\id\tens i\tens\id)(\id\tens R_\nabla)g\in \Omega^1 A\tens_A\Omega^1A\end{equation}
where $g\in  \Omega^1 A\tens_A\Omega^1 A$ and where we write its inverse as $\<,\>: \Omega^1 A\tens_A\Omega^1 A\to A$. We will check shortly that in the classical case this is -{1/2} of the usual Ricci tensor, i.e. the natural normalisation from a structural point of view is not quite the usual one. The Ricci scalar $S=\<\ ,\ \>({\rm Ricci})$ is likewise -{1/2} the usual one.

We now twist this. In view of the way that $\wedge^F$ is defined, clearly we can set
\begin{eqnarray*}
i^F(\xi\wedge\eta) &=& c\circ i(\xi\wedge\eta)=F^-(v_{[-1]}\tens w_{[-1]}) v_{[0]}\tens^F w_{[0]}
\end{eqnarray*}
if $i(\xi\wedge\eta)=v\tens w$ (or rather a sum of such terms) and descend to quotients. For example, if we start with classical geometry then clearly
\[ i^F(\xi\wedge\eta)={1\over 2}(F^-(\xi_{[-1]}\tens\eta_{[-1]})\,\xi_{[0]}\tens^F\eta_{[0]}-F^-(\eta_{[-1]}\tens\xi_{[-1]})\,\eta_{[0]}\tens^F \xi_{[0]}).\] 

We then compute the Ricci tensor defined analogously in the $\tens^F$ category,
 \begin{eqnarray*}
{\rm Ricci}^F=\Big((\<,\>^F\tens^F\id)\,\Phi^{-1}\,(\id\tens^F i^F)\tens^F\id\Big)\,
\Phi^{-1}\,(\id\tens^F R_{\nabla^F})\, g^F\ \in\ \Omega^1 A_F \tens_{A_F} \Omega^1 A_F\ .
\end{eqnarray*}

\begin{propos} \label{twistRicci}
${\rm Ricci}^F=c({\rm Ricci})$. The Ricci scalar  $S^F=S$ is unchanged under twist. 
\end{propos}
\noindent {\bf Proof:}\quad To compare the result with classical geometry,  we use the following commutative diagram, applying the top line to $g\to g^F$:
\begin{eqnarray}\label{vcusfghgakvcvuu}
\xymatrix{
\Omega^1 A \tens_A\Omega^1 A    \ar[rr]_{c}  \ar[d]_{\id\tens R_\nabla}       & &
    \Omega^1 A \tens^F_A\Omega^1 A \ar[d]_{\id\tens^F R_{\nabla^F}}      \\
    \Omega^1 A \tens_A(\Omega^2 A\tens_A\Omega^1 A)    \ar[rr]_{(\id\tens c)c}  \ar[d]_{=}       & &
    \Omega^1 A \tens^F_A(\Omega^2 A\tens^F_A\Omega^1 A) \ar[d]_{\Phi^{-1}}      \\
(\Omega^1 A \tens_A\Omega^2 A)\tens_A\Omega^1 A    \ar[rr]_{(c\tens \id)c}  \ar[d]_{(\id\tens i)\tens\id}       & &
(\Omega^1 A \tens^F_A\Omega^2 A)\tens^F_A\Omega^1 A \ar[d]_{(\id\tens^F i^F)\tens^F\id}      \\
(\Omega^1 A \tens_A(\Omega^1 A\tens_A\Omega^1 A))\tens_A\Omega^1 A  
  \ar[rr]_{((\id\tens c)c\tens\id)c}  \ar[d]_{=}       & &
(\Omega^1 A \tens^F_A(\Omega^1 A\tens^F_A\Omega^1 A))\tens^F_A\Omega^1 A \ar[d]_{\Phi^{-1}\tens^F\id}      \\
((\Omega^1 A \tens_A\Omega^1 A)\tens_A\Omega^1 A)\tens_A\Omega^1 A    \ar[rr]_{((c\tens \id)c\tens\id)c}  \ar[d]_{(\<,\>\tens\id)\tens\id}       & &
((\Omega^1 A \tens^F_A\Omega^1 A)\tens^F_A\Omega^1 A)\tens^F_A\Omega^1 A \ar[d]_{(\<,\>^F\tens^F\id)\tens^F\id}      \\
  \Omega^1 A \tens_A\Omega^1 A    \ar[rr]_{c}   & &
\Omega^1 A\tens^F_A\Omega^1 A    \\ 
}
\end{eqnarray}
and filling in the cells as commutative in view of the definitions of the various twisted objects in relation to the originals. The method here is the same as for other results above but at this point given with greater brevity.  \quad $\square$

This kind of argument is an explicit elaboration of the statement that the twist functor is an equivalence of categories. Nevertheless, there may be readers who would like an even more explicit treatment with tensors as familiar in Riemannian geometry. There are good reasons {\em not} to do this in noncommutative (and nonassociative) geometry, namely we would need first to elaborate a notion of {\em local coordinates}. We shall, however, check the definition of Ricci in (\ref{ricci}) in the classical case.  
We use standard formulae, for vector fields $U$, $V$ and $W$ and Christoffel symbols:
\begin{eqnarray*}
R(U,V)(W)^l &=&
R^l_{\phantom{l}ijk}\, W^i\, U^j\,V^k\ ,\cr
 R(U,V)(W)
 &=& \nabla_U\nabla_V\,W-\nabla_V\nabla_U\,W-\nabla_{[U,V]}\,W\ ,\cr
 (\nabla_V\,W)^l &=& V^k\,(W^l_{\phantom{l},k}+\Gamma^l_{ki}\,W^i)\ ,\cr
 R^{l}_{\phantom{l}ijk} &=& \frac{\del \Gamma^l_{ki}}{\del
x^j}\,-\, \frac{\del \Gamma^l_{ji}}{\del x^k}\,+\,
\Gamma^m_{ki}\,\Gamma^l_{jm}\,-\,\Gamma^m_{ji}\,\Gamma^l_{km}
\end{eqnarray*}
We use $\xi^a$ for the differential form $\extd x^a$ in terms of coordinates on
the manifold $M$. 
Then $\nabla\xi^a=-\Gamma_{bc}^a\,\xi^b\tens\xi^c$ and
\begin{eqnarray*}
R_\nabla(\xi^a) &=& (\extd\tens\id-\id\wedge\nabla)(-\Gamma_{bc}^a\,\xi^b\tens\xi^c) \cr
&=& -\, \Gamma_{bc,d}^a\,\xi^d\wedge \xi^b\tens\xi^c\ -\ 
\Gamma_{bc}^a\,\Gamma_{pq}^c\,\xi^b\wedge\xi^p\tens\xi^q  \cr
&=& -\, (\Gamma_{bc,d}^a + \Gamma_{de}^a\,\Gamma_{bc}^e)\,\xi^d\wedge \xi^b\tens\xi^c \cr
&=& -{1\over 2} R^a_{\phantom{a}cdb}\, \xi^d\wedge \xi^b\tens\xi^c\ .
\end{eqnarray*}
Using this to continue our calculation,
\begin{eqnarray*}
g &=&g_{fa}\,\xi^f \tens \xi^a\ ,\cr
(\id\tens R_\nabla)\, g &=& -{1\over 2}\,R^a_{\phantom{a}cdb}\,g_{fa}\,
\xi^f \tens
( \xi^d\wedge \xi^b\tens\xi^c)\ .
\end{eqnarray*}
When we apply $i$, we use the fact that we already have antisymmetry due to the
symmetries of $R^a_{\phantom{a}cdb}$, and so do not have to explicitly antisymmetrise.
Referring to (\ref{ricci}) we have
 \begin{eqnarray*}
{\rm Ricci} &=&
 -\, R^a_{\phantom{a}cdb}\,\,g_{fa}
 \,(\<,\>\tens\id)\,(\id\tens i)(\xi^f \tens
 \xi^d\wedge \xi^b)\tens\xi^c / 2\cr
 &=&
 -\, R^a_{\phantom{a}cdb}\,\,g_{fa}\, 
\,(\<,\>\tens\id)\,(\xi^f \tens
 (\xi^d\tens \xi^b))\tens\xi^c / 2\cr
   &=&
 -\, R^a_{\phantom{a}cdb}\,\,g_{fa}\, g^{fd}\,
 \xi^b\tens\xi^c / 2\cr
    &=&
 -\, R^a_{\phantom{a}cab}\,
 \xi^b\tens\xi^c/ 2\\ 
 &=& -{1\over 2}R_{cb}\xi^b\tens\xi^c\ 
\end{eqnarray*}
in our conventions (with $\tens_A$ understood throughout). The - sign has its roots in the fact that we work with Riemann as an operator on forms not on vector fields. The factor of {1/2} has its roots in the fact that $R_\nabla$ has values in $\Omega^2A$ rather than in antisymmetric tensors. Doing the same computations using the twisted objects similarly gives ${\rm Ricci}^F=-{1\over 2}R_{cb}F^-(\xi^b_{[-1]}\tens\xi^c_{[-1]})\xi^b_{[0]}\tens^F\xi^c_{[0]}$ as it should.

In view of the last proposition, we see that the Einstein tensor Ricci-${2\over n}g S$  (the trace-reversal of Ricci) also twists as ${\rm Einstein}^F=c({\rm Einstein})$. We can similarly extend the metric and $i$ to interior products, a Hodge *-operation,  etc, in particular to define the stress-energy tensor $T$ of matter fields of various types as a composition of maps, all of which twist. Then (in suitable normalisation)
\[ {\rm Einstein}= T\quad\Leftrightarrow\quad {\rm Einstein}^F=T^F.\]
Details will be elaborated elsewhere but the same remarks apply that this approach of itself works {\em too} well: we need experience and insight from noncommutative examples that go beyond twisting in order to know which of various ways to do these constructions is the most natural. 

\section{Some examples}

Examples of interest fall into two groups. The first consists of examples where the algebra of coordinates is associative and of independent interest. In this case knowing that such examples are twists opens the door to associated constructions for these algebra. Quite often in the literature, however, it is not realised that even when the algebra happens to be associative, this does not mean that the nonassociativity has gone away. Rather, as we explained in \cite{BM2}, it is merely {\em hidden} and typically resurfaces elsewhere, for example in the differential geometry.  The second class of examples are `manufactured' in the sense that the coordinate algebra, typically nonassociative, is defined for the first time, by twisting, and then studied. Since twisting works automatically, the content here is to have interesting classical data consisting of a geometry with a symmetry, which could be a discrete one, and an interesting cochain $F$.  We shall look in detail at examples of both types. Let us note that, although it is not our focus, we can also apply our theory to the case where $F$ is a cocycle and we stay within associative but noncommutative geometry, for example our results imply noncommutative conformal geometry on $\theta$-deformed spheres in the approach of \cite{BraMa}.

\subsection{Drinfeld twist examples}

We have already covered the first class of examples in \cite{BM1,BM2}, namely the coordinate Hopf algebras $C_q(G)$ associated to semisimple Lie algebras $\cg$,  using a different point of view on Drinfeld's  twist element $F$ that obtains this  as $C(G)^F$. We have explained the situation in the introduction. Unfortunately, the form of $F$ is not really known other than to lowest order. It indeed takes the form
\[ F(a\tens b)=\eps(a)\eps(b)+\lambda f(a\tens b)+O(\lambda^2)\]
where $f\in \cg\tens \cg$ and $\lambda$ is the deformation parameter. For matrix groups a Lie algebra element $\xi$ is evaluated against matrix element  functions $t^i{}_j\in C(G)$ by 
\[ \xi(t^i{}_j)=\rho^i{}_j(\xi)\]
where $\rho$ is a representation, and then extended to products via the additive coproduct on $\xi$. In Drinfeld's theory $F$ is given as a formal power-series of enveloping algebra elements but it can be similarly evaluated on elements to give formal power-series in $\lambda$. This was sufficient for the deformation-theory analysis in \cite{BM1} but it is expected that Drinfeld's $F$ defined in this way is actually finite and hence defined for numerical values of $\lambda$. The reason is that $F$ is quite similar in character to the `universal R-matrix' which makes $U_q(\cg)$ formally quasitriangular, and this does dualise to a well-defined linear functional making it coquasitriangular over the field. The underling reason is that the relevant elements in the power-series are nilpotent and hence some power of them vanish against any finite product of  representative functions. Moreover, from Drinfeld's theory one knows that $f_{21}-f=r_-$ the antisymmetric part of the Drinfeld-Sklyanin solution $r$ of the modified classical Yang-Baxter equations (which is known explicitly).  In this context it was shown in \cite{BM1} that while $C(G)^F$ happens to be associative, this necessarily does not apply to $\Omega(G)^F$. 

We note that although the $H^F$ construction means rather more, its algebra can be viewed as an example of the $(\ )_F$ algebra cochain twist. To do this, take $A=C(G)$ and $H=C(G)\tens C(G)^{cop}$ where left and  right translation are viewed as a left coaction of $H$. The standard $\Omega(G)$ is left-and-right translation invariant (bicovariant) and hence indeed twists by \[ F(a\tens b\tens c\tens d)=F(a\tens c) F^-(b\tens d)\]
where $F$ on the right is Drinfeld's $F$ viewed as a linear functional. We therefore have a twisted calculus $\Omega(G)_F$ for $C_q(G)$ as another point of view on \cite{BM1}.

In order to apply our further twisting theory we take the classical metric defined by the Killing form $K$, 
\[ g=K(\omega\tens_A\omega)\]
where $\omega\in \Omega^1(G)\tens \cg$ is the left Maurer-Cartan form and $K:\cg\tens\cg\to \C$ is supposed to be non-degenerate and ad-invariant. This is both left and right invariant and hence invariant under $H$. The Levi-Civita connection determined by this is
\[ \nabla \omega= -[\omega\tens_A\omega]\]
where we take the Lie bracket of the values of $\omega$. The latter obeys the Maurer-Cartan equations 
\[ \extd \omega+[\omega\wedge\omega]=0\]
and hence ${\rm Tor}_\nabla=0$. That $\nabla$ is metric compatible follows easily from ad-invariance of $K$ and a short computation yields $R_\nabla$ and ${\rm Ricci}\propto g$.  Now, because $K$ is ad-invariant the metric defined as above on terms of the left-invariant Maurer-Cartan form  is biinvariant, i.e. both left and right translations on the group are isometries. Hence $H=C(G)\tens C(G)^{cop}$ coacts and leaves the metric invariant and the results in Section~7 apply. Therefore we now also have a metric $g^F$ and Levi-Civita connection $\nabla^F$ etc. 

The above is an example of hidden non-associativity in that the algebra $A_F\isom C_q(G)$ is associative. 
We can also do the one-sided twist using only left-covariance of all the same structures and working directly with the dual version of Drinfeld's $F$ (not doubled up). Then as in \cite{BM2} we have $A_F=C(G)_F$ a non-associative algebra, with $\Omega(A_F)=\Omega(G)_F$, $g^F$ and $\nabla^F$ all defined by twisting. It remains, however, the case that working with Drinfeld's cochain is difficult due to lack of formulae.  

\subsection{Podles sphere and nonassociative black holes}

In this section we look at $\cg=sl_2$ and cochain $F$ of the form 
\[ F=1+\sum_{n=1}^\infty c_n(\lambda)( x\tens y)^n= f(x\tens y)\]
for some function $f$. Here $[x,y]=h$, $[h,x]=2x$, $[h,y]=-2y$ are $sl_2$ relations. Although for practical purposes one can think of $F\in U(su_2)\bar{\tens} U(su_2)$ we actually think of $F$ in the convolution algebra of linear functionals  $C(SU_2)\tens C(SU_2)\to \C$, or even better as an element of the smaller space $(C(SU_2)^*)^{\tens 2}$ i.e. as a tensor product of linear functionals.  Here $x$ (say) is regarded as such a functional and note that
\[ x^n(t^{i_1}{}_{j_1}\cdots t^{i_m}{}_{j_m})=0\]
for all $n>m$ due to $x^2(t^i{}_j)=0$ in the defining representation. With these remarks we are going to work formally with $F$ viewed in $U(su_2)\bar{\tens} U(su_2)$ but everything can be made precise in the above way. Clearly we have $(\eps\tens \id)F=(\id\tens\eps)F=1$ as needed for a cochain this form, and we require $F$ to be invertible. For this we require that the power-series $f$  has a non-zero radious of convergence about the origin. 

According to our theory any such cochain induces a noncommutative and non-associative geometry on any classical manifold which has an $SU_2$ symmetry. For Riemannian geometry it means any manifold $M$ with an $su_2$-triple of Killing vectors fields $\xi_i$ obeying,
\[ [\xi_i,\xi_j]=-\eps_{ijk}\xi_k\]
and we set $h\la  =2\imath\xi_3$, $x\la =-\xi_+$ and $y\la =-\xi_-$ where $\xi_\pm=\xi_1\pm\imath\xi_2$ to match our above conventions. Note that since we based our twists on left coactions this corresponds to a right action of $U(su_2)$ now. However, this Hopf algebra is isomorphic to its opposite so we can work with familiar left actions as stated. For our global theory we should suppose that these vector fields integrate to a right action of $SU_2$ and that everything is algebraic (so that we have a coaction). However, none of these is really necessary in deformation theory where we work only with power-series.

\begin{propos} We let $M=S^2$ defined as usual with coordinates $\sum x_i^2=1$ and $\xi_i=\eps_{ijk}x_j{\del\over\del x_k}$. Nontrivial $F$ of the form above can {\em never} induce an associative product. However, if $c_2=c_1^2$ the product of $C(S^2)_F$ coincides on the generators with that of the Podles sphere defined by generators $K,z,z^*$ with relations
\[ zK=q^2 Kz,\quad z^*z=1-K^2,\quad zz^*=1-q^2 K^2,\]
where $q^2=1-2c_1$. We assume that $c_1$ is real and $\ne 1, {1\over 2}$.
\end{propos}
\proof  We work in complex coordinates $x_\pm=x_1\pm\imath x_2$ and begin by computing  the action  as
\[ x\la\begin{pmatrix}x_+\cr x_-\cr  x_3\end{pmatrix}=\imath\begin{pmatrix}0\cr -2x_3\cr x_+\end{pmatrix},\quad
 y\la\begin{pmatrix}x_+\cr x_-\cr x_3\end{pmatrix}=\imath\begin{pmatrix}2x_3\cr 0\cr -x_-\end{pmatrix}.\]
 We compute $A_F$ from $a\bullet b= \cdot(F\la(a\tens b))=\cdot f(x\tens y)(a\tens b)$ to obtain
 \[ x_3\bullet x_3=x_3^2-c_1 x_+(-x_-),\quad x_+\bullet x_+=x_+^2,\quad x_-\bullet x_-=x_-^2\]
 \[ x_+\bullet x_-=x_+x_-,\quad  x_-\bullet x_+=x_-x_+-c_1(-2x_3)2 x_3-4c_2 x_+(-x_-)\]
 \[ x_+\bullet x_3=x_+x_3,\quad x_3\bullet x_+=x_3x_+-c_1x_+2x_3\]
 \[ x_-\bullet x_3=x_-x_3-c_1(-2x_3)(-x_-),\quad x_3\bullet x_-=x_3x_-.\]
 Now, from the last  four relations we see that $x_\pm$ indeed $q^2$-commute with $x_3$ as required for  $K\propto x_3$ and $z\propto x_-$. From the first relation and $x_+x_-=r^2-x_3^2$ for a sphere of radius $r$ we see that
 \[ x_3\bullet x_3=x_3^2+c_1(r^2-x_3^2)=(1-c_1)x_3^2+r^2c_1\] and hence $x_3^2(1-c_1)=x_3\bullet x_3-r^2 c_1$. We also see that
 \[ x_+\bullet x_-=x_+x_-=r^2-x_3^2={(r^2-x_3\bullet x_3)\over(1-c_1)}.\]
 Accordingly we set $z=x_-/\sqrt{1-c_1}$ and $K=x_3/r$ and have the 2nd relation of the Podles sphere. Note now that {\em if we assume associativity} among the generators,
 \[ x_-\bullet(x_+\bullet x_-)=x_-\bullet{(r^2-x_3\bullet x_3)\over 1-c_1}={x_-r^2\over 1-c_1}-{(1-2c_1)^2\over 1-c_1}(x_3\bullet x_3)\bullet x_-\]
 \[=(x_-\bullet x_+)\bullet x_-={(A-Bx_3\bullet x_3)\over (1-c_1)}\bullet x_-\]
 forces $A=r^2$ and $B=(1-2c_1)^2$ (in other words, the 3rd of the Podles sphere relations is dictated by the first two and associativity). Comparing with what we computed for $x_-\bullet x_+$ we need  $c_1^2=c_2$ and in this case we obtain 3rd relation as well.  A look at  $\bullet$ of some higher degree monomials, however, confirms that in this case  $\bullet$ is nevertheless not associative.\quad$\square$
 
We have included this example because it is striking that, while nonassociative, the nonassociativity is almost hidden and on the generators gives us a Podles sphere, i.e. a `mildly nonassociative' Podles sphere. Notice that the true Podles sphere is $C_q(SU_2)$-covariant which suggests that our form of $F$ has features similar to the Drinfeld cochain for $sl_2$. Also, the true Podles sphere has recently been interpreted as a $q$-fuzzy sphere\cite{Ma:fuz} or more geometrically as a slice in the $q$-hyperboloid of unit Lorentzian distance in Minkowski space and indeed arises as effective description of at least 3D quantum gravity {\em with} cosmological constant\cite{MaSch}.

We can apply this construction to spheres of any radius and to manifolds nested by spheres such as $\R^3$, the entire hyperboloid in Minkowski space, and the Schwarzschild black-hole. Clearly, it will induce a `mildly nonassociative' Podles or $q$-fuzzy sphere in place of each classical sphere. In particular, as the rotational symmetry is an isometry of the Schwarzschild solution, the metric is invariant and we have $g^F,\nabla^F$ on these spaces by our construction.

It is worth noting that these results are in the spirit of but rather different from recent attempts to apply {\em projective} Drinfeld twists to quantize black holes \cite{Schupp} based on a certain different $F=f(x\tens y)$ in \cite{Ale} with $c_n$ are generated by a beta-function. As we have seen, using this cochain in a standard way cannot give an associative algebra. Rather, the point of view is quite different; $F$ acts on $S^2=SU_2/U(1)$ via {\em right-translation} to induce a left-translation $SU_2$-invariant bidifferential operator. Here $U(su_2)$ itself does not act on $C(SU_2/U(1))$ in such a way at all, so this does not fit our framework. However, the failure of $F$ to be a cocycle is contained ($\phi=\del F$ has an $h$ to the right in its middle factor, and this acts trivially) and as a result the quantizaiton is accidentally associative, giving a fuzzy-sphere and a residual $SU_2$ symmetry. Although the framework is different from our more conventional use of Drinfeld twists, the lesson is the same that even though the algebra quantizing each sphere is associative, there remains a hidden nonassociativity in the model as $F$ is not a cocycle.
 
\subsection{Octonionic-$S^7$ as second quantisation}

In this section we give a first concrete example of nontrivial nonassociative Riemannian geometry, with curvature, of interest because it lies in the same category as the one in which the octonions live in the approach of \cite{AlbMa}. This is the category of $\Z_2^3$-graded spaces and we use the same cochain $F$ that `quantizes' $k(\Z_2^3)$ to obtain $\CO=k(\Z_2^3)_F$ as explained in Section~3. We again work  algebraically over a field $k$, and in the sequel we will assume this to be of characteristic zero (the reader can keep in mind $k=\R$). In physics when the result of first quantisation is (in some form) used as the basis of another (infinite dimensional) quantization, one speaks of {\em second quantization}. In the same spirit, we regard the result of the first (finite) quantization as a space on its own right and look at a submanifold of unit length octonions. This is  $S^7$  and we `quantize' it via a twist  albeit quite mildly by an action of $\Z_2^3$. These methods could be used to similarly twist any geometry admitting a finite group symmetry.

As usual, the effort is all in the preparation. We recall that the  description of octonions we use has basis $\{e_a\}$ labelled by $a\in\Z_2^3$.  The norm on the octonions in this basis is the usual euclidean one. Hence our natural description of $S^7$ is with coordinates $\{x^a\}$ again labelled  by $a\in \Z_2^3$ and the relation$\sum_a x_a^2=1$.  Following \cite{KliMa} we also have a coproduct, counit and antipode,
\[ \Delta x_c=\sum_{a+b=c} x_a \tens x_b F(a,b),\quad \eps x_a=\delta_{a,0},\quad S x_a=F(a,a)x_a\]
making $k(S^7)$ into a Hopf coquasigroup. It fails to be an ordinary Hopf algebra just because $S^7$ is not a group. Notice that all these maps and the algebra relations themselves respect a $\Z_2^3$ grading defined by $|x_a|=a$, i.e. $H=k\Z_2^3\isom k(\Z_2^3)$ coacts on everything here.

Next, in \cite{KliMa} Hopf-algebra-like methods are used to define the differential geometry of $S^7$ and we continue this programme now to compute its Riemannian geometry. The starting point is that a `standard' differential structure defined by ideal $(k(S^7)^+)^2$ as for algebraic groups, is parallelizable with basis of `left-invariant' 1-forms $\{\omega^i\}$ where $\omega^i=Sx^i\o\extd x^i\t=\sum_aF(a,i)x_a\extd x_{a+i}$ is defined as for Hopf algebras and computes as stated. Our convention is that middling indices $i,j,k,l,m,n$ run in the range $1,\cdots,7$. This induces `left-invariant vector fields' $\del_i$  and `Lie-like structure constants' $c^i{}_{jk}\in k(S^7)$  by
\[ \extd f=\sum_i (\del_i f)\omega^i,\quad \extd \omega^i+c^i{}_{jk}\omega^j\omega^k=0\]
which can be computed as
\[ \del_i=-\sum_a F(a,i) x^{a+i}{\del\over\del x^a}\]
\[ c^i{}_{jk}=-F(i,j)F(i+j,k)\sum_a\phi(a,i,j)\phi(a,i+j,k)F(a,i+j+k) x^a x^{a+i+j+k}.\]
Notice that the former do not close under commutator and the latter are not constants but functions, both effects due to the non-associativity of the product on $S^7$. Moreover,

\begin{propos} The `Killing form' is $c^m{}_{in} c^n{}_{jm}=-6\delta_{ij}$ (summation of repeated indices understood). Using $\delta_{ij}$ from now on to lower indices, $c_{ijk}$ is totally antisymmetric. \end{propos}
We refer to \cite{KliMa} for the methods used to obtain these results. However, $F$ and $\phi$ are given by particular formulae hence all our assertions can be checked  by direct computation using Mathematica. Being antisymmetric in the outer two indices then means that the metric
\[ g=\delta_{ij}\omega^i\tens\omega^j\]
is preserved by 
\[ \nabla \omega^i=-c^i{}_{jk}\omega^j\tens\omega^k\]
as 
\[ \nabla g=-c^i{}_{mn}\delta_{ij}\omega^m\tens\omega^n\tens\omega^j-c^j{}_{mn}\delta_{ij}\omega^i\tens\omega^m\tens\omega^n=0.\]
In calculations we leave $\tens_A$ understood. Here $\nabla$ acts as a derivation and its left output (which might be evaluated against a vector field) is kept to the far left and all $\tens$ products are over $k(S^7)$ so we can bring coefficients to the left as well. This connection is also manifestly torsion free and is therefore the Levi-Civita connection for the natural metric constructed in this way. These steps are exactly as for the Lie group case in Section~8.1 but now $c^i{}_{jk}$ are functions and do not define a Lie algebra (they define a Mal'tsev algebra at the identity of $S^7$), and we have given the computations in a basis rather than more abstractly. 

It remains to compute the Riemann and Ricci curvatures to make sure we have the usual constant curvature metric (or just for completeness). From (\ref{Rnabla}) we compute

 \begin{eqnarray*}
  R_\nabla\omega^i&=&(\extd\tens\id-\id\wedge\nabla)(-c^i{}_{jk}\omega^j\tens\omega^k)\\
  &=&-\extd c^i{}_{jk}\omega^j\tens\omega^k+c^i{}_{jk}c^j{}_{mn}\omega^m\omega^n\tens\omega^k-c^i{}_{jk}\omega^jc^k{}_{mn}\omega^m\tens\omega^n\\
  &=&-\left(c^i{}_{mj}c^j{}_{nk}-c^i{}_{jk}c^j{}_{mn}+\del_m c^i{}_{nk}\right)\omega^m\omega^n\tens\omega^k
  \end{eqnarray*}
 To obtain the Ricci curvature we apply a map $i$ that lifts the 2-form-valued part of $R_\nabla$ to an antisymmetrised element of $\Omega^1\tens\Omega^1$ (this is one way to be able to define an interior product). We then contract the first leg of the metric and the first leg of the lifting with the inverse metric $\<\omega^i,\omega^j\>=\delta_{ij}$ as in (\ref{ricci}),
 \begin{eqnarray*}
 {\rm Ricci}&=&\delta_{pi}\<\omega^p,(i\tens\id)R_\nabla\omega^i\>\\
 &=&-\delta_{pi}\<\omega^p,{1\over 2}\left(c^i{}_{mj}c^j{}_{nk}-c^i{}_{nj}c^j{}_{mk}-2c^i{}_{jk}c^j{}_{mn}+\del_m c^i{}_{nk}-\del_n c^i{}_{mk}\right)\omega^m\> \omega^n\tens\omega^k\\
 &=&-{1\over 2}\left(c^i{}_{ij}c^j_{nk}-c^i{}_{nj}c^j{}_{ik}-2c^i{}_{jk}c^j{}_{in}+\del_i c^i{}_{nk}-\del_n c^i{}_{ik}\right)\omega^n\tens\omega^k\\
 &=&{1\over 2}\left(c^i{}_{nj}c^j{}_{ki}-\del_ic^i{}_{nk}\right)\omega^n\tens\omega^k=-3g
 \end{eqnarray*}
 using the symmetries of $c$ and the Killiing form identity. Here the term $\del_i c^i{}_{nk}$ is antisymmetric in $n,k$ hence it must vanish as the Ricci tensor is necessarily symmetric. This turns out to be the case as one can compute using Mathematica from the formulae already given. Hence this $S^7$ is an Einstein space as it should be (the usual answer here for a unit 7-sphere would be $6g$; the minus sign and factor of 1/2 is due to our conventions as explained at the end of Section~7). In short, the main difference that the $c^i{}_{jk}$ are functions does show up in the full Riemann tensor but  not in Ricci. 
 
 We have given an algebraic derivation of the usual geometry of $S^7$ using Hopf coquasigroup methods. This is actually a uniform approach that applies to all the parallellizable spheres. Now notice that $|\omega^i|=i$ (here $\extd$ does not change the $\Z_2^3$-degree) and $|c^i{}_{jk}|=i+j+k$. Hence $g$ has zero degree and $\nabla$ does change the degree, i.e. everything is covariant under the grading Hopf algebra $k\Z_2^3$ and we can apply our cochain twist. We obtain $k(S^7)_F$ and a calculus $\Omega(S^7)_F$ on it, with everything nonassociative. The product $x^a\bullet x^b=F(a,b)x^ax^b$ obeys
 \[ (x^a\bullet x^b)\bullet x^c=\phi(a,b,c) x^a\bullet(x^b\bullet x^c),\quad x^b\bullet x^a=\CR(a,b)x^a\bullet x^b\]
 where $\phi(a,b,c)=-1$ {if and only if} $a,b,c$ are linearly independent, as before and $\CR(a,b)=-1$ {if and only if} $a,b$ are linearly independent, and otherwise $1$, again as for the octonions. The second feature means that the algebra is {\em altercommutative} in the sense that elements anticommute {if and only if} they have distinct non-zero degrees. Otherwise they commute. This is not, however, the octonion algebra. For products of generators we have degree (and the same for their $\bullet$ products)
 \[ |(x^0)^{m}(x^1)^{n_1}\cdots (x^7)^{n_7}|=d(N)=(n_4+n_5+n_6+n_7,n_2+n_3+n_6+n_7,n_1+n_3+n_5+n_7)\]
where $N=(n_1,\cdots,n_7)$ and the degree is computed modulo 2 to give an element of $\Z_2^3$. We have enumerated the $\{x^a\}$ in binary order $x^{000},x^{001},\cdots,x^{111}$ for brevity. Then
 \[( \prod (x^i)^{n_i})\bullet (\prod (x^j)^{m_j})=F(d(N),d(M)) \prod (x^i)^{n_i+m_i}\]
 and associativity remains via $\phi$ on the relevant degrees (so three monomials associate {if and only if} their $\Z_2^3$-degrees are linearly independent) and commutativity remains via $\CR$ (so two monomials commute {if and only if} they do not have different non-zero degrees). The generator $x^0$ has zero degree, has unchanged product with everything and any power of it associates with everything. Finally, we have the sphere relations which now takes the form
 \[ 1=x^0\bullet x^0-\sum_i x^i\bullet x^i.\]
The even powers $x^i\bullet x^i=-(x^i)^2$ and $1,x^0$ generate a commutative associative subalgebra of degree 0, i.e. which commutes and associates with everything else. Similarly, we can compute the differential calculus obtained by twisting the classical one. The products
\[ x^a \bullet \extd x^b=F(a,b)x^a\extd x^b,\quad \extd x^a \bullet \extd x^b=F(a,b)\extd x^a\extd x^b\]
obey parallel formulae controlled by $\phi,\CR$. Explicitly,
\[ 0=x^0\bullet\extd x^0-\sum_ix^i\bullet\extd x^i\]
with relations involving $x^0$ or $\extd x^0$ unchanged, as is $\extd x^a\bullet\extd x^a=0$ for all $a$. The changed relations for distinct $i,j,k$ are
\[ x^i\bullet \extd x^j=-\extd x^j\bullet x^i,\quad \extd x^i\bullet \extd x^j=\extd x^j\bullet \extd x^i,\quad (x^i\bullet\extd x^j)\bullet x^k=\phi(i,j,k)x^i\bullet (\extd x^j\bullet\extd x^k)\]
and similar. In degree three we have
\[ (\extd x^i\bullet\extd x^j)\bullet \extd x^k=-\extd x^i\bullet(\extd x^j\bullet\extd x^k)\]
{if and only if} $i+j+k\ne 0$. 

By our theory the metric and connection also twist:
\[ g^F=\omega^0\tens\omega^0-\sum_i\omega^i\tens\omega^i,\quad \nabla^F\omega^i=\sum_{j,k}F(i,k)c^i{}_{jk}\omega^i\tens\omega^k\]
where we made use of $F(i+k,k)=F(i,k)F(k,k)=-F(i,k)$ to simplify. We write the summations explicitly for clarity. The inverse metric $\<\ , \>^F$ also acquires $-1$ values on the $i=1,\cdots,7$ diagonal entries and 
\[ i^F(\omega^i\bullet\omega^j)={1\over 2}(\omega^i\tens\omega^j+\omega^j\tens\omega^i),\quad i\ne j\]
because $F(i,j)=-F(j,i)$ when we compute $i^F(\omega^i\omega^j)$, which we can then cancel as the same factor occurs in $\omega^i\bullet\omega^j$. Following through the calculation of ${\rm Ricci}_F$ in terms of these quantities and remembering the associator, one can verify
\[ {\rm Ricci}^F=-3 g^F.\]
as implied by  Propostion~\ref{twistRicci} (we again have an Einstein space, with constant curvature). This provides a useful check of our nonassociative theory.

Further insight into this `octonionic' $S^7$ can be gained by looking at the same construction for second-quantized quaternions and complex numbers. The quaternions are constructed in the same way but with $\Z_2^2$-grading. $F$ has the same form but without the `cubic' term. As a result, $F$ is a cocycle and $\phi=1$. The above geometry gives the classical Riemannian geometry of $S^3$ viewed as the unit quaternions, and $k(S^3)_F$ is associative and altercommutative. It therefore has bosonic central generators $1,x^0$, non-commutative generators $x^1,x^2,x^3$ with relations
\[ \{x^1,x^2\}=\{x^2,x^3\}=\{x^3,x^1\}=0, \quad 1=(x^0)^2-(x^1)^2-(x^2)^2-(x^3)^2\]
where the bullet product is understood. The noncommutative generators are `Grassman-like' but do not anticommute with themselves, indeed the algebra is commutative in the category of $\Z_2^2$-graded spaces with a specific symmetric braiding different from that of the category of supervector spaces.  If we go down one level more to $\Z_2$-gradings we have one degree 1 generator $x^1$ and $1,x^0$ are bosonic. As a result everything is commutative and associative and $k(S^1)_F=k[x^0,x^1]$ modulo the relation $1=(x^0)^2-(x^1)^2$, i.e.  a hyperbola. Again one has metric $g^F$ and connection $\nabla^F$ with constant curvature. So all of these `twisted spheres' are actually some kind of hyperbolae. One could view them as baby versions of the $q$-deformation twist in Section~8.2 which also turned spheres into mildly nonassociative $q$-hyperboloids.

\end{document}